\documentclass[10pt]{amsart}

\usepackage{mathrsfs}
\usepackage{enumerate}
\usepackage{color}
\usepackage{relsize}
  \usepackage{ifpdf}
  \ifpdf
    \usepackage[all]{xy}
    \SelectTips{cm}{}
    \usepackage[pdftex]{graphicx}
    \usepackage[pdftex,pagebackref,colorlinks,linkcolor=blue,citecolor=blue,urlcolor=blue,a4paper,hypertexnames=true,plainpages=false]{hyperref}
    \hypersetup{pdfauthor={Clara L\"oh and Roman Sauer}, pdftitle={Degree theorems and Lipschitz simplicial volume for non-positively curved manifolds of finite volume}}
  \else
    \usepackage[all, dvips]{xy}
    \SelectTips{cm}{}
    \usepackage{graphicx}
    \usepackage[pagebackref,colorlinks,linkcolor=blue,citecolor=blue,urlcolor=blue,a4paper,hypertexnames=true,plainpages=false]{hyperref}
  \fi

\usepackage{amsrefs}

\newcommand{\bbN}{{\mathbb N}}
\newcommand{\bbQ}{{\mathbb Q}} 
\newcommand{\bbR}{{\mathbb R}} 
\newcommand{\bbH}{{\mathbb H}} 
\newcommand{\bbZ}{{\mathbb Z}} 

\newcommand{\bfG}{\mathbf{G}}

\newcommand{\calu}{\mathcal{U}}
\newcommand{\calr}{\mathcal{R}}

\newcommand{\scrC}{\mathscr{C}}
\newcommand{\scrH}{\mathscr{H}}

\let\epsilon\varepsilon
\let\rho\varrho
\let\theta\vartheta
\let\phi\varphi

\newcommand{\abs}[1]{{\left\lvert #1\right\rvert}}
\newcommand{\absfix}[1]{{\bigl\lvert #1\bigr\rvert}}
\newcommand{\supn}[1]{\|#1\|_\infty}
\newcommand{\bs}{\backslash}
\newcommand{\cd}{\operatorname{cd}}
\newcommand{\colim}{\operatorname*{colim}}
\newcommand{\const}{\operatorname{const}}

\newcommand{\defq}{\mathrel{\mathop:}=}
\newcommand{\dvol}{\operatorname{dvol}}
\newcommand{\id}{\operatorname{id}}
\newcommand{\im}{\operatorname{im}}

\newcommand{\isom}{\operatorname{Isom}}

\newcommand{\lfvol}[1]{\norm{#1}}
\newcommand{\lipschitz}{\operatorname{Lip}}
\newcommand{\lip}{\mathrm{Lip}}
\newcommand{\lipvol}[1]{\norm{#1}_{\mathrm{Lip}}}
\newcommand{\lipminvol}{\operatorname{minvol_\lip}}
\newcommand{\simvol}[1]{\|#1\|}
\newcommand{\map}{\operatorname{map}}
\newcommand{\minvol}{\operatorname{minvol}}
\newcommand{\norm}[1]{{\left\lVert #1\right\rVert}}

\newcommand{\ricci}{\operatorname{Ricci}}
\newcommand{\rk}{\operatorname{rk}}
\newcommand{\sd}{\operatorname{sd}}
\newcommand{\smear}{\operatorname{smear}}
\newcommand{\supp}{\operatorname{supp}}
\newcommand{\str}{\operatorname{str}}
\newcommand{\vol}{\operatorname{vol}}
\DeclareMathOperator{\diam}{diam}
\newcommand{\catzero}{{\upshape{\smaller CAT}}$(0)$}
\newcommand{\frface}[2]{#2\rfloor_{#1}} 
\newcommand{\baface}[2]{{}_{#1}\lfloor#2}
\newcommand{\fa}[1]{\forall_{#1}\spc} 
\newcommand{\spc}{\;\;\;}
\newcommand{\ucov}[1]{\widetilde{#1}}
\newcommand{\args}{\,\cdot\,}
\newcommand{\invlim}{\varprojlim}
\newcommand{\EZ}{\mathrm{EZ}}
\newcommand{\AW}{\mathrm{AW}}
\DeclareMathOperator{\SL}{SL}
\DeclareMathOperator{\SO}{SO}

\DeclareMathOperator*{\esssup}{ess-sup}

\renewcommand{\colon}{:}

\newcommand{\Hde}{H_{\mathrm{dR}}}                    
\newcommand{\Clonelip}{C^{\smash{\ell^1},\mathrm{Lip}}}

\newcommand{\Clf}{C^{\mathrm{lf}}}           
\newcommand{\Hlf}{H^{\mathrm{lf}}}           
\newcommand{\Clfl}{C^{\mathrm{lf},<L}}       
\newcommand{\Hlfl}{H^{\mathrm{lf},<L}}       
               
\newcommand{\Ccs}{C_{\mathrm{cs}}}           
\newcommand{\Hcs}{H_{\mathrm{cs}}}           
\newcommand{\Ccsl}{C_{\mathrm{cs},<L}}       
\newcommand{\Hcsl}{H_{\mathrm{cs},<L}}       
\newcommand{\Clflip}{C^{\mathrm{lf},\mathrm{Lip}}}             
\newcommand{\Hlflip}{H^{\mathrm{lf},\mathrm{Lip}}}             
\newcommand{\Ccslip}{C_{\mathrm{cs},\mathrm{Lip}}}             
\newcommand{\Hcslip}{H_{\mathrm{cs},\mathrm{Lip}}}             
\newcommand{\Clip}{C^{\mathrm{Lip}}}         
\newcommand{\Hlip}{H^{\mathrm{Lip}}}         

\newcommand{\Cmealip}{\scrC^{\mathrm{Lip}}}        
\newcommand{\Hmealip}{\scrH^{\mathrm{Lip}}}        
\newcommand{\lfset}{S^{\mathrm{lf}}}     
\newcommand{\liplfset}{S^{\mathrm{lf}, \mathrm{Lip}}}     
\DeclareMathOperator{\Str}{Str}              

\newcommand{\smallcoprod}{\sqcup}

\newtheorem{theorem}{Theorem}[section]
\newtheorem{lemma}[theorem]{Lemma}
\newtheorem{proposition}[theorem]{Proposition}
\newtheorem{corollary}[theorem]{Corollary}
\theoremstyle{definition}
\newtheorem{definition}[theorem]{Definition}

\newtheorem{example}[theorem]{Example}
\newtheorem{remark}[theorem]{Remark}

\newtheorem{setup}[theorem]{Setup}
\numberwithin{equation}{section}

 \makeatletter
  \let\c@equation\c@theorem
  \renewcommand{\tagform@}[1]{\maketag@@@{#1}}
  \let\c@figure\c@theorem
\makeatother

\title[Lipschitz simplicial volume of non-positively curved
  manifolds]{Degree theorems and Lipschitz simplicial volume for 
  non-positively curved manifolds\\ of finite volume}

\author{Clara L\"oh}
\address{Graduiertenkolleg ``Analytische Topologie und
  Metageometrie,''  
  Westf\"alische Wil\-helms-Universit\"at M\"unster,
  M\"unster, Germany} 
\email{clara.loeh@uni-muenster.de}
\urladdr{wwwmath.uni-muenster.de/u/clara.loeh}

\author{Roman Sauer}
\address{University of Chicago, Chicago, USA}
\curraddr{Westf\"alische Wilhelms-Universit\"at M\"unster, M\"unster, Germany}
\email{romansauer@member.ams.org}
\urladdr{www.romansauer.de}

\date{\today}

\begin{document}

\begin{abstract}
  We study a metric version of the simplicial volume on Riemannian
  manifolds, the Lipschitz simplicial volume, with applications to
  degree theorems in mind. We establish a proportionality 
  principle and a product inequality from which we derive an 
  extension of Gromov's volume comparison theorem to 
  products of negatively curved manifolds or locally symmetric 
  spaces of non-compact type. In contrast, we provide vanishing 
  results for the ordinary simplicial volume; for instance, we show
  that the ordinary simplicial volume of non-compact locally symmetric
  spaces with finite volume of $\bbQ$-rank at least~$3$ is zero. 
\end{abstract}

\maketitle

\section{Introduction and statement of results}\label{sec:intro}

\noindent
The prototypical \emph{degree theorem} bounds the 
degree~$\deg f$ of a proper, continuous map $f:N\rightarrow M$ 
between $n$-dimensional Riemannian manifolds of finite volume by 
\[\deg(f)\le \const_n \cdot \frac{\vol(N)}{\vol(M)}. \]
 
For example, Gromov's volume comparison theorem~\cite{gromov}*{p.~13}
is a degree theorem where the target~$M$ has negative sectional
curvature and the domain~$N$ satisfies a lower Ricci curvature bound.
In \emph{loc.~cit.} Gromov also pioneered the use of the
\emph{simplicial volume} to prove theorems of this kind. Recall that
the simplicial volume~$\norm{M}$ of a manifold $M$ without boundary is
defined by
\[  \norm{M}=
    \inf\bigl\{\abs{c}_1;~\text{$c$ fundamental cycle of $M$ with $\bbR$-coefficients}
        \bigr\}.
\]
Here $\abs{c}_1$ denotes the $\ell^1$-norm with respect to 
the basis given by the singular simplices. If $M$ is non-compact 
then one takes locally finite fundamental cycles in the 
above definition. Under the given curvature assumptions, Gromov's comparison 
theorem is proved by the following three steps (of which the third  
one is elementary): 
\begin{enumerate}
\item Upper volume estimate for target: $\norm{M}\ge\const_n\vol(M)$.
\item Lower volume estimate for domain: $\norm{N}\le\const_n\vol(N)$. 
\item Degree estimate: $\deg(f)\le \const_n\norm N/\norm M$. 
\end{enumerate} 

\emph{Unless stated otherwise, all manifolds in this text are 
assumed to be connected and without boundary. 
As Riemannian metrics on locally symmetric spaces of non-compact type
we always choose the standard metric, i.e., the one given by the Killing
form~\cite{eberlein}*{Section~2.3.11}.}

\subsection{Main results}\label{subsec:main results}

In this article, we prove degree theorems where the target is
non-positively curved and has finite volume. More specifically, we
consider the case where the target is a product of negatively curved
manifolds of finite volume or locally symmetric spaces of finite
volume. To this end, we study a variant of the simplicial volume, the
\emph{Lipschitz simplicial volume}, and pursue a Lipschitz version of
the three step strategy above. The properties of the Lipschitz
simplicial volume we show en route are also of independent interest. 

Before introducing the Lipschitz simplicial volume, we give a brief
overview of the properties of the ordinary simplicial volume of
non-compact locally symmetric spaces of finite volume: On the one hand
by a classic result of Thurston~\cite{thurston}*{Chapter~6} the
simplicial volume of finite volume hyperbolic manifolds is
proportional to the Riemannian volume. According to Gromov and
Thurston the simplicial volume of complete Riemannian manifolds with
pinched negative curvature and finite volume is
positive~\cite{gromov}*{Section~0.3}.  In addition, we proved by
different means that the simplicial volume of Hilbert modular
varieties is positive~\cite{spinoff} (see also
Theorem~\ref{thm:hilbert modular varieties} below).  In accordance
with these examples we expect positivity for all locally symmetric
spaces of $\bbQ$-rank~$1$.

On the other hand, in Section~\ref{sec:vanishing results} we show that
the simplicial volume of locally symmetric spaces of $\bbQ$-rank at
least~$3$ vanishes -- in particular, the ordinary simplicial volume
does not give rise to the desired degree theorems.

\begin{theorem}\label{thm:vanishing of locally finite simplicial volume}
Let $\Gamma$ be a torsion-free, arithmetic lattice of a semi-simple, center-free  
\mbox{$\bbQ$-group}~$\bfG$ with no compact factors. 
Let $X=\bfG(\bbR)/K$ be the associated symmetric 
space where $K$ is a maximal compact subgroup of~$\bfG(\bbR)$. 
If $\Gamma$ has \mbox{$\bbQ$-rank} at least~$3$, then $\lfvol{\Gamma\bs X}=0$. 
\end{theorem} 

This result is based on a more general vanishing theorem
(Corollary~\ref{cor: small classifying space}) derived from Gromov's
vanishing-finiteness theorem~\cite{gromov}*{Corollary~(A) on p.~58} by
constructing suitable amenable coverings for manifolds with nice
boundary and whose fundamental groups admit small classifying spaces.

Gromov's original applications of the vanishing-finiteness theorem
contain the surprising fact that the simplicial volume of any product
of three open manifolds is zero~\cite{gromov}*{p.~59}. However, there
are products of two open manifolds whose simplicial volume is non-zero
(see Example~\ref{ex:products with non-zero sv}), and Gromov's
argument fails for products of two open surfaces. In particular, the
$\bbQ$-rank~$2$ case is still open.

In contrast to Theorem~\ref{thm:vanishing of locally finite simplicial
  volume}, Lafont and Schmidt showed the following positivity
result in the closed case~\cite{ben+jean}; the proof is based on work of
Connell-Farb~\cite{chris+benson}, as well as -- for the exceptional
cases -- Thurston, Savage, and Bucher-Karlsson:

\begin{theorem}[Lafont, Schmidt]\label{thm:simplicial volume of closed locally
    symmetric spaces}
  Let $M$ be a closed locally symmetric space of non-compact
  type. Then $\norm M > 0$.
\end{theorem}

In view of the fact that the simplicial volume of non-compact
manifolds is zero in a large number of cases, Gromov studied 
geometric variants of the simplicial
volume~\cite{gromov}*{Section~4.4f}, i.e., simplicial volumes where
the simplices allowed in fundamental cycles respect a geometric
condition. In this article, we consider the following Lipschitz
version of simplicial volume: 

\begin{definition}\label{def:metric simplicial volume}
  Let $M$ be an $n$-dimensional, oriented Riemannian manifold. 
  For a locally finite chain $c\in \Clf_n(M)$ we denote
  the supremum of the Lipschitz constants of the simplices occurring
  in~$c$ by $\lipschitz(c) \in [0, \infty]$.  The \emph{Lipschitz simplicial
  volume}~$\lipvol{M}\in [0,\infty]$ of~$M$ is defined by
  \begin{equation*}
          \lipvol{M}
    = \inf\bigl\{ \abs{c}_1
               ;~\text{%
		      $c \in \Clf_n(M)$ fundamental cycle of~$M$ with
                      $\lipschitz(c) < \infty$}
              \bigr\}. 
  \end{equation*}
\end{definition}

By definition, we have the obvious inequality $\lfvol{M} \le
\lipvol{M}$. It is easy to see that 
if $f:N\rightarrow M$ is a proper Lipschitz map 
between Riemannian manifolds, then 
\[\deg(f)\cdot\lipvol{M}\le\lipvol{N}.\]

\begin{remark}\label{rem:lip equals nonlip}
  If $M$ is a \emph{closed} Riemannian manifold, then
  $\norm{M}=\lipvol{M}$; each fundamental cycle involves only finitely
  many simplices, and hence this equality is implied by the fact that
  singular homology and smooth singular homology are isometrically
  isomorphic~\cite{loeh}*{Proposition~5.3}.
\end{remark}

In Section~\ref{sec:proportionality for non-compact manifolds} 
we prove the following theorem, which leads to  
a degree theorem for locally symmetric spaces of finite volume.  

\begin{theorem}[Proportionality principle]
  \label{thm:proportionality most general}
  Let $M$ and $N$ be complete, non-positively curved Riemannian
  manifolds of finite volume. Assume that their
  universal covers are isometric. Then
  \[   \frac{\lipvol{M}}{\vol(M)}
     = \frac{\lipvol{N}}{\vol(N)}.
  \]
\end{theorem}

The proportionality principle for closed Riemannian manifolds is a
classical theorem of Gromov~\citelist{\cite{gromov}*{Section~0.4}
  \cite{thurston}*{pp.~6.6--6.10}\cite{loehdiplom}*{Chapter~5}}.
The proportionality principle in the closed case does not
require a curvature condition, and our proof in the non-closed case 
uses non-positive curvature in a light way.  
It might be possible to weaken the curvature condition in the
non-compact case. 

By Theorem~\ref{thm:simplicial volume of closed locally symmetric
  spaces} the proportionality principle for the ordinary simplicial
volume cannot hold in general since for every locally symmetric space
of finite volume there is always a \emph{compact} one such that their
universal covers are isometric~\cite{borelcompact}.  For the same
reason, Theorems~\ref{thm:proportionality most general}
and~\ref{thm:simplicial volume of closed locally symmetric spaces} and
Remark~\ref{rem:lip equals nonlip} imply the following corollary.

\begin{corollary}\label{cor:lipschitz simplical volume positiv}
The Lipschitz simplicial volume of locally symmetric spaces of finite 
volume and non-compact type is non-zero. 
\end{corollary}

Gromov~\cite{gromov}*{Section~4.5} states also a proportionality
principle for non-compact manifolds for geometric invariants related
to the Lipschitz simplicial volume. Unraveling his definitions, one
sees that it implies a proportionality principle for finite volume
manifolds without a curvature assumption (which we need) provided one
of the manifolds is compact (which we do not need). This would be
sufficient for the previous corollary.  Gromov's proof, which is
unfortunately not very detailed, and ours seem to be independent.

The simplicial volume of a product of oriented, closed, connected
manifolds can be estimated from above as well as from below in terms
of the simplicial volume of both
factors~\citelist{\cite{gromov}*{p.~17f}\cite{bp}*{Theorem~F.2.5}}. While
the upper bound continues to hold for the locally finite simplicial
volume in the case of non-compact
manifolds~\cite{loehphd}*{Theorem~C.7}, the lower bound in general
does not.

The Lipschitz simplicial volume on the other hand is better behaved 
with respect to products. In addition to the estimate~$\lipvol{M \times N}
\leq c(\dim M + \dim N) \cdot \lipvol M \cdot \lipvol N$, the presence
of non-positive curvature enables us to derive also the non-trivial
lower bound:

\begin{theorem}[Product inequality for non-positively curved manifolds]
 \label{thm:product formula}
  Let $M$ and $N$ be two complete, non-positively curved Riemannian
  manifolds. Then
  \[ \lipvol M \cdot \lipvol N
     \leq \lipvol{M \times N}.
  \] 
\end{theorem}

On a technical level, we mention two issues that often prevent one
from extending properties of the simplicial volume for compact
manifolds to non-compact ones, and thus force one to work with the
Lipschitz simplicial volume instead.  Firstly, there is no
straightening (see Section~\ref{subsec:geodesic straightening}) for
locally finite chains: The straightening of a locally finite chain $c$
is not necessarily locally finite. However, it is locally finite
provided $\lip(c)<\infty$, which motivates a Lipschitz condition.
Secondly, there is no well-defined cup product for compactly supported
cochains.  This is an issue arising in the proof of the product
inequality.  We circumvent this difficulty by introducing the complex of
cochains with Lipschitz compact support (see
Definition~\ref{def:lipschitz dual}), which carries a natural
cup-product.

\subsection{Degree theorems}\label{subsec:applications}

To apply the theorems of the previous section to degree theorems, 
we need upper and lower estimates of the volume by the Lipschitz 
simplicial volume. 

For the (locally finite) simplicial volume and all complete 
$n$-dimensional Riemannian manifolds, 
Gromov gives the bound $\norm{M}\le (n-1)^nn!\vol(M)$ provided 
$\ricci(M)\ge -(n-1)$~\cite{gromov}. The latter stands for 
$\ricci(M)(v,v)\ge -(n-1)\norm{v}^2$ for all~$v \in TM$. 
One can extract from \emph{loc.~cit.} a similar estimate 
for the Lipschitz simplicial volume: 

\begin{theorem}[Gromov]\label{thm: volume and Lipschitz simplicial volume}
For every $n\ge 1$ there is a constant $C_n>0$ such that 
every complete $n$-dimensional Riemannian manifold $M$ 
with sectional curvature~$\sec(M)\le 1$ and Ricci curvature 
$\ricci(M)\ge -(n-1)$ satisfies 
\[ \lipvol{M}\le C_n\cdot\vol(M).\]
\end{theorem}

\begin{proof}
  For $\sec(M)\le 0$ this follows from~\cite{gromov}*{Theorem (A),(4),
    in Section~4.3} by applying it to~$U=M$, $R=1$, a fundamental
  cycle $c$, and $\epsilon\rightarrow 0$: One obtains a fundamental
  cycle $c'$ made out of straight simplices whose diameter is less
  than~$R+\epsilon$.  In particular, $\lip(c')<\infty$ by
  Proposition~\ref{prop:geodesic simplex - non-positive curvature}.
  Further, the estimate $\lipvol{M}\le\norm{c'}\le C_n\vol(M)$ follows
  from (4) in \emph{loc.~cit.} and the Bishop-Gromov inequality, which
  provides a bound of~$l_v'(R)$ in terms of~$n$~\citelist{\cite{gromov}*{(C) in
  Section~4.3}\cite{gallot}*{Theorem~4.19 on~p.~214}}.

Gromov also explains why these arguments carry over to the general
case that $\sec(M)\le 1$~\cite{gromov}*{Remarks~(B) and (C) in
  Section~4.3}. In this case, $c'$ is made out of straight simplices of
diameter less than~$\pi/2$ (Section~\ref{subsec:geodesic simplices}), and
$\lip(c')<\infty$ follows from Proposition~\ref{prop:geodesic simplex
  - bounded from above}.
\end{proof}

\begin{corollary}\label{cor:finite lipschitz volume}
Any complete Riemannian manifold of finite volume that has 
an upper sectional curvature and lower Ricci curvature bound 
has finite Lipschitz simplicial volume. 
\end{corollary}

Connell and Farb~\cite{chris+benson} prove, building upon techniques
of Besson-Courtois-Gallot, a degree theorem where the target $M$ is a
locally symmetric space (closed or finite volume) with no local
$\bbR$, $\bbH^2$, or $\SL(3,\bbR)/\SO(3,\bbR)$-factor.  For
non-compact $M$ they have to assume that $f:N\rightarrow M$ is
(coarse) Lipschitz. Using the simplicial volume (and the work by
Connell-Farb, Thurston, Savage, and Bucher-Karlsson), Lafont and
Schmidt~\cite{ben+jean} prove degree theorems for closed locally symmetric
spaces including the exceptional cases. The following theorem
includes also the non-compact exceptional cases.

\begin{theorem}[Degree theorem, complementing~\cites{chris+benson, ben+jean}]
\label{thm:degree theorem for locally symmetric spaces}
For every $n\in\bbN$ there is a constant $C_n>0$ with the following
property: Let $M$ be an $n$-dimensional locally symmetric space of
non-compact type with finite volume. Let $N$ be an $n$-dimensional
complete Riemannian manifold of finite volume with $\ricci(N)\ge
-(n-1)$ and $\sec(N)\le 1$, and let $f:N\rightarrow M$ be a proper
Lipschitz map. Then
\[\deg(f)\le C_n\cdot\frac{\vol(N)}{\vol(M)}.\]
\end{theorem}

\begin{proof}
  By Theorem~\ref{thm:proportionality most general} and
  Corollary~\ref{cor:lipschitz simplical volume positiv} we know that
  $\lipvol{M}=\const_n\vol(M)$ where $\const_n>0$ depends only on the
  symmetric space~$\smash{\ucov{M}}$.  Because there are only finitely
  many symmetric spaces (with the standard metric) in each dimension,
  there is $D_n>0$ depending only on $n$ such that $\lipvol{M}\ge
  D_n\vol(M)$. So Theorem~\ref{thm: volume and Lipschitz simplicial
    volume} applied to~$N$ and $\lipvol{N}\ge\deg(f)\lipvol{M}$
  yield the assertion.
\end{proof}

Unfortunately, the Lipschitz simplicial volume cannot be used to prove
positivity of Gromov's \emph{minimal volume} $\minvol(M)$ of a smooth
manifold~$M$; the minimal volume is defined as the infimum of volumes
$\vol(M,g)$ over all complete Riemannian metrics~$g$ on~$M$ whose
sectional curvature is pinched between~$-1$ and~$1$.

Next we describe the appropriate modification of $\minvol(M)$ in our
setting: The \emph{Lipschitz class} $[g]$ of a complete
Riemannian metric~$g$ on~$M$ is defined as the set of all complete
Riemannian metrics~$g'$ such that the identity $\id: (M,g')\rightarrow
(M,g)$ is Lipschitz. Then we define the \emph{minimal volume of $[g]$}
as
\[\lipminvol(M,[g])=\bigl\{\vol(M,g')
  ;~\text{$-1\le \sec(g')\le 1$ and  $g'\in [g]$}\bigr\}.\]
Of course, we have $\lipminvol(M,[g])=\minvol(M)$ whenever  
$M$ is compact. 
Theorem~\ref{thm:degree theorem for locally symmetric spaces},  
applied to the identity map and varying metrics, implies: 

\begin{theorem}
\label{thm:minimal volume of locally symmetric spaces}
The minimal volume of the Lipschitz class of the standard  
metric of a locally symmetric space of non-compact type and finite
volume is positive. 
\end{theorem}

Excluding certain local factors, Connell and Farb 
have the following stronger statement for the minimal volume 
instead of the Lipschitz minimal volume. 

\begin{theorem}[Connell-Farb]
The minimal volume of a locally symmetric space of non-compact type and 
finite volume 
that has no local $\bbH^2$- or
$\mathrm{SL}(3,\bbR)/\mathrm{SO}(3,\bbR)$-factors is positive. 
\end{theorem}

A little caveat: Connell and Farb state this theorem 
erroneously as a corollary of a degree theorem for which they have to assume a
Lipschitz condition. This would only give the positivity of the Lipschitz 
minimal volume. However, Chris Connell explained to us how 
to modify their proof to get the positivity of the minimal volume. 

As an application of the product inequality we obtain a new 
degree theorem for products of manifolds 
with (variable) negative curvature or locally symmetric spaces. 

\begin{theorem}[Degree theorem for products]
\label{thm:degree theorem for products}
For every~$n\in\bbN$ there is a constant~$C_n>0$ with the following
property: Let $M$ be a Riemannian $n$-manifold of finite volume that
decomposes as a product~$M = M_1 \times \dots \times M_m$ of
Riemannian manifolds, where for every~$i\in\{1,\ldots,m\}$ the
manifold $M_i$ is either negatively curved with $-\infty<-k <
\sec(M_i)\le -1$ or a locally symmetric space of non-compact type.
Let $N$ be an $n$-dimensional, complete Riemannian manifold of finite
volume with~$\sec(N)\le 1$ and $\ricci(N)\ge -(n-1)$. Then for every
proper Lipschitz map~$f:N\rightarrow M$ we have
\[\deg(f)\le C_n\cdot\frac{\vol(N)}{\vol(M)}.\]
\end{theorem}

\begin{proof}
In the sequel, $D_i,D_i', E_n$, and $C_n$ stand for 
constants depending only on~$n$. 
If $M_i$ is negatively curved then 
Thurston's
theorem~\citelist{\cite{gromov}*{Section~0.3}\cite{thurston}} 
yields 
\[\vol(M_i)\le D_n\lfvol{M_i}\le D_n\lipvol{M_i}.\]
If $M_i$ is locally symmetric of non-compact type 
then, as in the proof of 
Theorem~\ref{thm:degree theorem for locally symmetric spaces}, we 
also obtain $\vol(M_i)\le D'_n\lipvol{M_i}$. 
By the product inequality (Theorem~\ref{thm:product formula}),  
\[\vol(M)\le \max_{i\in\{1,\dots,m\}}(D_i,D_i')^m\lipvol{M}.\]
On the other hand, by Theorem~\ref{thm: volume and Lipschitz
simplicial volume}, we have $\lipvol{N}\le E_n\vol(N)$.  
Combining everything with~$\lipvol{N}\ge\deg(f)\lipvol{M}$, 
proves the theorem with the constant~$C_n=E_n/\max_{i\in\{1,\dots,m\}}(D_i,D_i')^m$. 
\end{proof}

As a concluding remark, we mention a computational application of the
proportionality principle.  We proved that $\norm{M}=\lipvol{M}$ for
Hilbert modular varieties~\cite{spinoff}. This fact combined with the
proportionality principle~\ref{thm:proportionality most general} and
work of Bucher-Karlsson~\cite{michelle} leads then to the following
computation~\cite{spinoff}:

\begin{theorem}\label{thm:hilbert modular varieties}
Let $\Sigma$ be a non-singular Hilbert modular surface. Then 
\[\norm{\Sigma}=\frac{3}{2\pi^2}\vol(\Sigma). \]
\end{theorem}

Conversely, the proportionality principle~\ref{thm:proportionality
  most general} together with Thurston's computation of the simplicial
volume of hyperbolic manifolds shows that the simplicial volume of
hyperbolic manifolds of finite volume equals the Lipschitz simplicial
volume. More generally, this holds true for locally 
symmetric spaces of 
$\bbR$-rank $1$~\cite{spinoff}*{Theorem~1.5; see also beginning of Section~1.5}. 
However, in the general $\bbQ$-rank~$1$ case, the relation between
the simplicial volume and the Lipschitz simplicial volume remains open. 

\subsection*{Organization of this work}

Section~\ref{sec:straightening} reviews the basic properties of
geodesic simplices and Thurston's straightening. The product inequality
(Theorem~\ref{thm:product formula}) is proved in
Section~\ref{sec:product formula}. Section~\ref{sec:proportionality
  for non-compact manifolds} contains the proof of the proportionality
principle (Theorem~\ref{thm:proportionality most general}). Finally,
Section~\ref{sec:vanishing results} is devoted to the proof of the
vanishing result (Theorem~\ref{thm:vanishing of locally finite
  simplicial volume}).

\subsection*{Acknowledgements}
The first author would like to thank the Graduiertenkolleg
``Analytische Topologie und Metageometrie'' at the WWU M\"unster for
its financial support.  The second 
author acknowledges support of the German Science Foundation
(DFG), made through grant SA 1661/1-1.

Both authors thank the University of Chicago for 
an enjoyable working atmosphere (C.L. visited UC in March/April~2007). 
We are very grateful to Chris Connell and Benson Farb for discussions 
about their work~\cite{chris+benson}. 
The second author also thanks Juan Souto and 
Shmuel Weinberger for several helpful discussions.  

\setcounter{tocdepth}{1}
\setcounter{tocdepth}{2}

\section{Straightening and Lipschitz estimates of straight simplices}\label{sec:straightening}

\noindent 
In Section~\ref{subsec:geodesic simplices}, we 
collect some basic properties of geodesic simplices. 
We recall the technique of straightening singular chains for 
non-positively curved manifolds in 
Section~\ref{subsec:geodesic straightening}. Variations of this
straightening play an important role in the proofs of the proportionality
principle (Theorem~\ref{thm:proportionality most general}) and the
product inequality (Theorem~\ref{thm:product formula}).

\subsection{Geodesic simplices}\label{subsec:geodesic simplices}

Let $M$ be a simply connected, complete Riemannian manifold. 
Firstly assume that $M$ has non-positive sectional curvature. 
For points~$x$,~$x'$ in $M$, we denote 
by~$[x,x'] \colon [0,1]\rightarrow M$ the unique geodesic joining~$x$ and~$x'$. 
The \emph{geodesic join} of two maps~$f$ 
and~$g \colon X \rightarrow M$ from a space $X$ to $M$ 
is the map defined by   
\[ [f,g] \colon X\times [0,1]\rightarrow M, (x,t)\mapsto\bigl[ f(x), g(x)\bigr](t).
\]
We recall the notion of geodesic simplex: 
The standard simplex $\Delta^n$ is given by 
$\Delta^n=\{(z_0,\ldots,z_n)\in\bbR^{n+1}_{\ge 0};~\sum_i z_i=1\}$, and 
we identify $\Delta^{n-1}$ with the subset~$\{(z_0,\ldots,z_n)\in\Delta^n;~z_n=0\}$. 
Moreover, the standard simplex is always equipped with the induced 
Euclidean metric. 
Let $x_0, \dots, x_n \in M$. 
The \emph{geodesic simplex}~$[x_0, \dots, x_n] \colon \Delta^n\rightarrow M$ 
with vertices~$x_0,\dots, x_n$ is defined inductively as 
\[   [x_0,\ldots, x_n]\bigl((1-t)s+t(0,\ldots,0,1)\bigr) 
   = \bigl[[x_0, \dots, x_{n-1}](s), x_n\bigr](t)
\]
for~$s\in\Delta^{n-1}$ and~$t \in [0,1]$.

More generally, if $M$ admits an upper bound~$K_0 \in (0, \infty)$ of
the sectional curvature, then every pair of points with distance less
than~$K_0^{-1/2}\pi/2$ in~$M$ is joined by a unique geodesic. Thus we
can define the \emph{geodesic simplex} with vertices $x_0,\dots, x_n$
as before whenever $\{x_0,\dots,x_n\}$ has diameter less
than~$K_0^{-1/2}\pi/2$~\cite{gromov}*{4.3~(B)}.

In the following two sections (Sections~\ref{subsubsec:estimates for
joins} and~\ref{subsubsec:lipschitz constant}), we provide uniform
estimates for Lipschitz constants of geodesic joins and simplices. 

\subsubsection{Lipschitz estimates for geodesic
  joins}\label{subsubsec:estimates for joins}

\begin{proposition}\label{prop:smooth join}
  Let $M$ be a simply connected, complete Riemannian manifold of
  non-positive sectional curvature, and let $n \in\bbN$. Let $f$,~$g
  \in \map(\Delta^n,M)$ be smooth maps. Then $[f,g]$ is smooth and has
  a Lipschitz constant that depends only on the Lipschitz constants
  for $f$ and $g$.
\end{proposition}

\begin{proof}
Using the exponential map we can 
rewrite $[f,g]$ as 
\begin{equation}\label{eq:smoothness of geodesic join}
[f,g](x,t)=\exp_{f(x)}\bigl(t\cdot\exp_{f(x)}^{-1}(g(x))\bigr).  
\end{equation}
Since the exponential map viewed as a map $TM\rightarrow M\times M$ 
is a diffeomorphism, $[f,g]$ is smooth. 
The assertion about the Lipschitz constant 
is a consequence of the following lemma. 
\end{proof}

\begin{lemma}\label{lem:cat0 estimate}
  Let $X$ be a compact metric space and $M$ as above. 
If $f$ and~$g\colon X \rightarrow M$ are two
  Lipschitz maps, then the geodesic join~$[f,g] \colon X\times [0,1]
  \rightarrow M$ is also a Lipschitz map, and we have
  \[ \lipschitz [f,g]
     \leq 2 \cdot 
          \bigl( \lipschitz f + \lipschitz g 
	       + \diam (\im f \cup \im g)
	  \bigr).
  \]
\end{lemma}

\begin{proof}
  Let $(x,t)$,~$(x',t') \in X \times [0,1]$. The triangle
  inequality yields
  \begin{align*}
              d_M\bigl( [f,g](x,t), [f,g](x',t')
                 \bigr)
       \leq &\; d_M\bigl( [f(x),g(x)](t), [f(x),g(x)](t') 
		 \bigr)\\
          + &\; d_M\bigl( [f(x),g(x)](t'), [f(x'),g(x')](t')
	         \bigr). 
  \end{align*}
  Because $\lipschitz[f(x),g(x)] =  d_M(f(x),g(x))$, the first term satisfies
  \begin{align*}
              d_M\bigl( [f(x),g(x)](t), [f(x),g(x)](t') 
		 \bigr)
       \leq &\; |t-t'| \cdot d_M\bigl( f(x), g(x)
                              \bigr)\\
       \leq &\; |t-t'| \cdot \diam\bigl(\im f \cup \im g
                                \bigr);
  \end{align*}
  notice that $\diam(\im f \cup \im g)$ is finite because $X$ is
  compact. 
  The \catzero-inequality allows us to simplify the second term as
  follows
  \begin{align*}
              d_M\bigl( [f(x),g(x)](t'), [f(x'),g(x')](t')
	         \bigr)
       \leq &\; d_M\bigl( [f(x),g(x)](t'), [f(x),g(x')](t')
		 \bigr)\\
	  + &\; d_M\bigl( [f(x),g(x')](t'), [f(x'),g(x')](t')
	         \bigr)\\
       \leq &\; (1-t') \cdot d_M\bigl( g(x), g(x')\bigr)\\
          + &\; t'     \cdot d_M\bigl( f(x), f(x')\bigr)\\
       \leq &\; d_M\bigl( g(x), g(x')\bigr)
          +   d_M\bigl( f(x), f(x')\bigr)\\
       \leq &\; \lipschitz f \cdot d_X(x,x')
            + \lipschitz g \cdot d_X(x,x').
  \end{align*}
  Therefore, we obtain
  \begin{align*}
              d_M\bigl( [f,g](x,t), [f,g](x',t')
                 \bigr)
     \leq  &\;      \bigl( \lipschitz f + \lipschitz g 
	                 + \diam (\im f \cup \im g)
	            \bigr)\\
     \cdot &\;     2 \cdot 
                   d_{X\times [0,1]} 
                    \bigl( (x,t), (x',t')
                    \bigr)
              . \qedhere
  \end{align*}
\end{proof}

\subsubsection{Lipschitz estimates for geodesic simplices}\label{subsubsec:lipschitz constant}

Similarly to geodesic joins also geodesic simplices admit a uniform
Lipschitz estimate and analogous smoothness properties. 

\begin{proposition}\label{prop:geodesic simplex - non-positive
    curvature}
Let $M$ be a complete, simply connected, non-positively curved Riemannian
manifold. Then every geodesic simplex in~$M$ is smooth. Moreover, for
every $D>0$ and $k\in\bbN$ there is $L>0$ such that every geodesic
$k$-simplex $\sigma$ of diameter less than~$D$ satisfies
$\norm{T_x\sigma}<L$ for every $x\in\Delta^k$.
\end{proposition}

\begin{remark}\label{rem:diameter of geodesic simplices}
  Let $M$ be a simply connected, complete Riemannian manifold of
  non-positive sectional curvature. If $x_0 , \dots, x_k \in M$, then
  applying the triangle inequality inductively shows that 
  \[ \fa{y \in \Delta^k} 
     d_M\bigl([x_0, \dots, x_k](y), x_k\bigr)
     \leq k \cdot \max_{i,~j\in\{0,\dots,k\}} d_M(x_i, x_j)
  \]
  and hence that
  \[      \diam\bigl(\im [x_0, \dots, x_k]\bigr)
     \leq 2 \cdot k \cdot \max_{i,~j\in\{0,\dots,k\}} d_M(x_i, x_j).
  \]
\end{remark}

In the proof of Theorem~\ref{thm: volume and Lipschitz simplicial
 volume}, it is necessary to have a more general version of 
Proposition~\ref{prop:geodesic simplex - non-positive curvature} 
dealing with a positive upper sectional curvature bound. 
In this case, locally, the same arguments apply: 

\begin{proposition}\label{prop:geodesic simplex - bounded from above}
Let $M$ be a complete, simply connected Riemannian manifold whose sectional 
curvature is bounded from above by~$K_0 \in (0,\infty)$. Then 
every geodesic
simplex $\sigma$ of diameter less than~$K_0^{-1/2}\pi/2$ is smooth. 
Further, there is a constant~$L>0$ such that 
every geodesic
$k$-simplex $\sigma$ of diameter less than~$K_0^{-1/2}\pi/2$ 
satisfies $\norm{T_x\sigma}<L$ for every $x\in\Delta^k$. 
\end{proposition}

The proofs of the following two lemmas used to prove
Propositions~\ref{prop:geodesic simplex - non-positive curvature}
and~\ref{prop:geodesic simplex - bounded from above} are elementary
and thus omitted. The proof of the first one is very similar to Lee's
proof of the Sturm comparison theorem~\cite{leerie}*{Proof of
  Theorem~11.1}.

\begin{lemma}\label{lem:differential equation}
Let $u:[0,1]\rightarrow\bbR_{\ge 0}$ be a smooth 
function 
such that $u(0)=0$, and $u(t)>0$ for $t\in (0,1]$, as well as
\[\fa{t \in [0,1]}\frac{d^2}{dt^2}u(t)+ \frac{\pi^2}4 \cdot u(t)\ge 0.\]
Then for all $t\in [0,1]$ we have 
\[u(t)\le u(1)\cdot \sin\bigl(t \cdot \pi/2\bigr).\] 
\end{lemma}

\begin{lemma}\label{lem:linear algebra}
Let $f:V\rightarrow W$ be a linear map between finite-dimensional 
vector spaces with inner products. Let $H\subset V$ be a subspace 
of co-dimension $1$, and let $z\in V$ be a vector such that 
$z$ and $H$ span $V$. Let $\{y_1,\ldots,y_{k-1}\}$ be an 
orthonormal basis of $H$. 
Assume that for some~$C > 0$
\[ \fa{w\in\{z,y_1,\ldots,y_{k-1}\}} 
   \norm{f(w)}\le C\cdot\norm{w}.
\]
Further, assume that the angle $\alpha$ between $z$ and $H$ lies in 
$[\epsilon, \pi/2]$ with $0<\epsilon\le\pi/2$. 
Then there is a constant $L>0$ that depends only on 
$\dim(V)$, $C$, and $\epsilon$ such that 
\[\norm{f}< L.\]
\end{lemma}

\begin{proof}[Proof of 
Proposition~\ref{prop:geodesic simplex - non-positive curvature} 
and~\ref{prop:geodesic simplex - bounded from above}]
That geodesic simplices are smooth is easily seen using the 
fact that the exponential map is a diffeomorphism. 
Let $K_0\ge 0$ be an upper bound for the sectional curvature 
of~$M$. By normalizing the metric we may assume that either $K_0=0$ 
or $K_0=1$. In the case~$K_0=1$, it is understood that~$D=\pi/2$. 
Led by the inductive definition of geodesic simplices, we prove the
proposition by induction over~$k$: For $k = 0$ or $k = 1$ there is
nothing to show. 

We now assume that there is an~$L' > 0$ such that every geodesic
$(k-1)$-simplex of diameter less than~$D$ is smooth and that the norm
of its differential is less than~$L'$.  Let $\sigma
:=[x_0,\ldots,x_k]:\Delta^k\rightarrow M$ be a geodesic $k$-simplex of
diameter less than~$D$. By 
the induction hypothesis, 
\begin{equation}\label{eq:induction hypothesis}
  \fa{p \in \Delta^{k-1}}
  \bigl\|T_p[x_0,\ldots,x_{k-1}]\bigr\| < L'.
\end{equation} 

In the following, we write~$v_0, \dots, v_k$ for the vertices
of~$\Delta^k$. Let $p\in\Delta^{k-1}$, and let $\gamma: [0,1]\rightarrow
M$ denote the geodesic from $x_k$ to $[x_0, \dots, x_{k-1}](p)$. 
Choose an orthonormal basis~$\{X_1,\ldots, X_{k-1}\}$ of the hyperplane
in~$\bbR^{k}$ spanned by~$\Delta^{k-1}$; then we can view 
$\{X_1,\ldots,X_{k-1}\}$ as an orthonormal frame of $T\Delta^{k-1}$.  
For $i\in\{1,\ldots, k-1\}$ we consider the following variation of
$\gamma$: 
\begin{equation*}
  H_i:(-\epsilon_i,\epsilon_i)\times [0,1]\rightarrow M,
  ~(s,t) \mapsto \bigl[ x_k, \sigma(s\cdot X_i + p) \bigr](t). 
\end{equation*}
Let $X_i(t):=\frac{d}{ds} H_i(s,t)\vert_{s=0}\in TM$. By definition, 
$X_i$ is a Jacobi field along~$\gamma$. Moreover, we have at
each point~$p(t):= [v_k, p](t)$ of~$\Delta^k$ the relation
\begin{equation}\label{eq:image under differential}
 T_{p(t)}\sigma(X_i)=t \cdot X_i(t).
\end{equation}
In order to obtain the desired bounds for~$\|T_{p(t)}\sigma\|$ we
first give estimates for~$\|X_i(t)\|$ and then apply
Lemma~\ref{lem:linear algebra} to conclude the proof.

For the following computation, let $D_t$ denote the covariant
derivative along $\gamma$ at~$\gamma(t)$, 
and let $K$ and $R$ denote the
sectional curvature and the curvature tensor, respectively. 
Straightforward differentiation and 
the Jacobi equation yield 
\begin{align*}
       \frac{d^2}{dt^2}\bigl\|X_i(t)\bigr\|^2 
  &=   2\cdot \norm{D_tX_i (t)}^2 
     - \Bigl\langle R\Bigl(X_i (t),\frac{d}{dt}\gamma\Bigr)
       \frac{d}{dt}\gamma, X_i (t)\Bigr\rangle\\
  &\ge -K_0
        \cdot \norm{X_i (t)}^2
        \cdot \Bigl\|\frac{d}{dt} \gamma\Bigr\|^2\\
  &\ge -K_0
        \cdot \norm{X_i (t)}^2
        \cdot D^2. 
\end{align*}
By definition, $X_i(0) = 0$, and 
by~\eqref{eq:image under differential} and~\eqref{eq:induction hypothesis}, 
\[\norm{X_i(1)}=\norm{T_p\sigma(X_i)}=
\norm{T_p[x_0,\ldots,x_{k-1}](X_i)} < L'.\] 
First assume that $K_0=0$. Then the smooth function $t\mapsto \norm{X_i(t)}^2$ 
starts with the value $0$, is non-negative, and convex. So it 
is non-decreasing. This implies that 
\begin{equation}\label{eq:estimate of Jacobi field}
  \fa{t \in [0,1]}
  \norm{X_i(t)}\le \norm{X_i(1)} < L'.
\end{equation} 
Next assume that $K_0=1$, thus $D=\pi/2$. 
Lemma~\ref{lem:differential equation} 
yields~\eqref{eq:estimate of Jacobi field}. 
Thus, in both cases $K_0=0$ or $K_0=1$ 
we see that $\norm{X_i(t)}\le L'$ for all $t\in
[0,1]$.  Further note that
\[     \Bigl\|T_{p(t)}\sigma\Bigl(\frac{d}{dt}p\Bigr)\Bigr\|
  =    \Bigl\|\frac{d}{dt}\gamma\Bigl\|
  \leq D.
\]
So Lemma~\ref{lem:linear algebra} implies 
that there is a constant $L>0$ that depends only on $L'$, $D$, and 
$k$ such that 
\[\|T_{\gamma(t)}\sigma\| < L \]
because the angle between the line $p(t)$ and $\Delta^{k-1}$ 
is at least $\epsilon>0$ with 
$\epsilon$ depending only on $\Delta^k$. 
\end{proof}

\subsection{Geodesic straightening}\label{subsec:geodesic straightening}

In the following, we recall the definition of the geodesic straightening 
map on the level of chain complexes, as introduced by 
Thurston~\cite{thurston}*{p.~6.2f}.

Let $M$ be a connected, complete Riemannian manifold of non-positive
sectional curvature.  A singular simplex on~$M$ is \emph{straight} if
it is of the form~$p_M\circ \sigma$ for some geodesic simplex~$\sigma$
on~$\ucov M$, where $p_M \colon \ucov M \rightarrow M$ is the
universal covering map.  The subcomplex of the singular
complex~$C_*(M)$ generated by the straight simplices is denoted
by~$\Str_*(M)$; the elements of~$\Str_*(M)$ are called \emph{straight
  chains}.  Every straight simplex is uniquely determined by the
(ordered set of) vertices of its lift to the universal cover.

The \emph{straightening}~$s_M \colon
  C_*(M) \rightarrow \Str_*(M)$ is defined by
  \[    s_M(\sigma) 
     :=  p_M \circ 
	 \bigl[ \ucov\sigma(v_0), \dots, \ucov\sigma(v_*)
	 \bigr]\text{ for $\sigma\in\map(\Delta^*,M)$},
  \]
  where $p_M\colon \ucov M \rightarrow M$ is the universal covering
  map, $v_0, \dots, v_*$ are the vertices of~$\Delta^*$, and $\ucov
  \sigma$ is some $p_M$-lift of~$\sigma$.

Notice that the definition of~$s_M(\sigma)$ is independent of the
chosen lift~$\ucov \sigma$ because the fundamental group~$\pi_1(M)$
acts isometrically on~$\ucov M$.

\begin{proposition}[Thurston]\label{prop:geodesic
  straightening} 
  Let $M$ be a connected, complete Riemannian manifold of non-positive
  sectional curvature.  Then the straightening~$s_M \colon C_*(M)
  \rightarrow \Str_*(M)$ and the inclusion~$\Str_*(M) \rightarrow
  C_*(M)$ are mutually inverse chain homotopy equivalences.
\end{proposition}

The easy proof is based on Lemma~\ref{lem:lemma from lee} below, which
is a standard device for constructing chain homotopies. Because we
need this lemma later, we reproduce the short argument for
Proposition~\ref{prop:geodesic straightening} here.

\begin{proof}[Proof of Proposition~\ref{prop:geodesic straightening}]
  For each singular simplex~$\sigma \colon \Delta^n \rightarrow M$
  on~$M$, we define
  \[ H_\sigma := p_M \circ
                 \bigl[ \ucov\sigma
                      , [\ucov\sigma(v_0), \dots, \ucov\sigma(v_n)]
                 \bigr]
     \colon \Delta^n \times [0,1] \rightarrow M,
  \]
  where $v_0, \dots, v_n$ are the vertices of~$\Delta^n$, and $\ucov
  \sigma$ is a lift of~$\sigma$ with respect to the universal covering
  map~$p_M$. It is not
  difficult to see that $H_\sigma$ is independent of the chosen
  lift~$\ucov\sigma$ and that $H_\sigma$ satisfies the hypotheses of
  Lemma~\ref{lem:lemma from lee} below. 

  Therefore, Lemma~\ref{lem:lemma from lee} provides us with a chain
  homotopy between~$\id_{C_*(M)}$ and the straightening map~$s_M$. 
\end{proof}

\begin{lemma}\label{lem:lemma from lee}
Let $X$ be a topological space. For each $i\in\bbN$ and each singular
$i$-simplex $\sigma:\Delta^i\rightarrow X$ let
$H_\sigma:\Delta^i\times I\rightarrow X$ be a homotopy such that 
for each face map $\partial_k:\Delta^{i-1}\rightarrow\Delta^i$ we have 
\begin{equation*}
H_{\sigma\circ\partial_k}=H_\sigma\circ (\partial_k\times\id_I).
\end{equation*}
Then $f^{(0)}$ and $f^{(1)}:C_\ast(X)\rightarrow C_\ast(X)$, defined by
$f^{(m)}(\sigma)=H_\sigma\circ i_m$ for~$m \in \{0,1\}$, are chain maps.
For every $i\in\bbN$ there are $i+1$ affine simplices 
$G_{k,i}:\Delta^{i+1}\rightarrow\Delta^i\times I$ such that 
\begin{equation*}
H:C_i(X)\rightarrow C_{i+1}(X),~h(\sigma)=\sum_{k=0}^iH_\sigma\circ
G_{k,i}
\end{equation*}
defines a chain homotopy $f^{(0)}\simeq f^{(1)}$. 
\end{lemma}

\begin{proof}
  This is literally proved in Lee's book~\cite{lee}*{Proof of
  Theorem~16.6, \mbox{p.~422-424}} although the lemma above is not stated
  as such.
\end{proof}

\begin{remark}\label{rem:notation from lee}
  The simplices $G_{k,i}$ in the previous lemma arise from decomposing
  the prism~$\Delta^i\times I$ into $(i+1)$-simplices.
\end{remark}

\section{Product inequality for the Lipschitz simplicial
  volume}\label{sec:product formula}

\noindent
This section is devoted to the proof of the product inequality  
(Theorem~\ref{thm:product formula}). 

The corresponding statement in the compact case is proved by first
showing that the simplicial volume can be computed in terms of bounded
cohomology and then exploiting the fact that the cohomological
cross-product is compatible with the semi-norm on bounded
cohomology~\citelist{\cite{gromov}*{p.~17f}\cite{bp}*{Theorem~F.2.5}}.
In a similar fashion, the product inequality for the locally finite
simplicial volume can be shown if one of the factors is
compact~\citelist{\cite{gromov}*{p.~17f}\cite{loehphd}*{Appendix~C}}.

To prove the Lipschitz version, we proceed in the following steps:
\begin{enumerate}
  \item We show that the Lipschitz simplicial volume can be computed
        in terms of a suitable semi-norm on cohomology with Lipschitz
        compact supports; this semi-norm is a variant of the supremum
        norm parametrized by locally finite supports
        (Sections~\ref{subsec:lflip homology},~\ref{subsec:cslip
        cohomology}, and~\ref{subsec:lipvol via cohomology}).
  \item The failure of the product inequality for the locally finite simplicial 
        volume is linked to the fact that there is no well-defined cross 
        product on compactly supported cochains. In contrast, we show in  
        Lemma~\ref{lem:cross product restricts} in 
        Section~\ref{subsec:cross-products} that there is a cross-product 
        for cochains with Lipschitz compact support 
        (Definition~\ref{def:lipschitz dual}), and 
        we analyze the interaction between this semi-norm and the
        cross-product on cohomology with compact supports
        (Section~\ref{subsec:cross-products}).
  \item Finally, we prove that the presence of non-positive curvature
        allows us to restrict attention to locally finite fundamental
        cycles of the product that have nice supports
        (Section~\ref{subsec:sparse}). This enables us to use the
        information on cohomology with Lipschitz compact supports to
        derive the product inequality (Section~\ref{subsec:product
        formula conclusion}).
\end{enumerate}

\subsection{Locally finite homology with a Lipschitz constraint}
\label{subsec:lflip homology}
The locally finite simplicial volume is defined in terms of the
locally finite chain complex. In the same way, the Lipschitz
simplicial volume is related to the chain complex of chains with
\emph{Lipschitz locally finite support}. 

\begin{definition}\label{def:family of compact subsets}
For a topological space~$X$, we define~$K(X)$ to be the set of all 
compact, connected, non-empty subsets of~$X$. 
\end{definition}

For simplicity, we consider only connected compact subsets. This is
essential when considering relative fundamental classes of pairs of
type~$(M, M-K)$.

\begin{definition}\label{def:locally finite lipschitz homology}
  Let $X$ be a metric space, and let $k \in \bbN$. Then we write
  \begin{align*}     
         \lfset_k (X)
     &:= \bigl\{ A \subset \map(\Delta^k,X)
               ;~\fa{K \in K(X)}
                |\{ \sigma \in A \mid \im(\sigma) \cap K \neq \emptyset \}| < \infty 
         \bigr\}\\
         \liplfset_k(X)
     &:= \bigl\{ A\in\lfset_k(X)
               ;~\exists_{L \in\bbR_{>0}}\spc
                 \fa{\sigma\in A}
                 \lipschitz(\sigma) < L
         \bigr\}.
  \end{align*}

  The elements of~$\liplfset_k(X)$ are said to be \emph{Lipschitz locally
  finite}. The subcomplex of~$\Clf_*(X)$ of all chains with Lipschitz
  locally finite support is denoted by~$\Clflip_*(M)$, and the
  corresponding homology -- so-called \emph{homology with Lipschitz
  locally finite support} -- is denoted by~$\Hlflip_*(X)$.
\end{definition}

\begin{theorem}\label{thm:locally finite lipschitz homology gives the locally finite homology} 
  Let $M$ be a connected Riemannian manifold. Then the homomorphism
  \[ \Hlflip_*(M)
     \rightarrow
     \Hlf_*(M)
  \]
  induced by the inclusion~$\Clflip_*(M) \rightarrow \Clf_*(M)$ is an
  isomorphism.  
\end{theorem}

During the course of the proof of this theorem, we rely on the
following notation:

\begin{definition}
  Let $X$ be a proper metric space, and 
  let $A \subset X$ be a subspace. Let $L \in\bbR_{>0}$. 
  \begin{enumerate}
    \item We write~$\Clfl_*(X)$ for the subcomplex of~$\Clflip_*(X)$
          given by
          \[    \Clfl_\ast(X)
	     := \bigl\{ c \in \Clf_\ast(X)
	              ;~\lipschitz(c) < L
		\bigr\}.
          \]
    \item Similarly, we define $C_*^{<L}(X) := \{ c \in
          C_*(X);~\lipschitz(c) < L\}$ as well as $C_*^{<L}(X,A) :=
          C_*^{<L}(X)/C_*^{<L}(A)$.
    \item We use the abbreviation
          \[ \Clip_*(X) := \colim_{L \rightarrow \infty} 
                           C_*^{<L}(X)
                        =  \bigl\{ c\in C_*(X)
			         ;~\lipschitz(c) < \infty
			   \bigr\}.
	  \]
    \item The corresponding homology groups are denoted
          by~$\Hlfl_*(X)$,~$H_*^{<L}(X)$, $H_*^{<L}(X,A)$,
          and~$\Hlip_*(X)$ respectively.
  \end{enumerate}
\end{definition}

By definition, we can express the chain complex of chains with locally
finite Lipschitz support via the colimit 
\[  \Clflip_*(X)\xleftarrow{\cong} 
    \colim_{L \rightarrow \infty} \;\Clfl_*(X).
\]
with the obvious inclusions as structure maps. 
Moreover, if $X$ is connected, the term on the right hand side
expands to the inverse limit 
\[  \fa{L \in \bbR_{>0}}
    \Clfl_*(X)\xrightarrow{\cong}
    \invlim_{K \in K(X)} C_*^{<L}(X, X - K)
\]
with the obvious projections as structure maps. 
\begin{proof}[Proof of Theorem~\ref{thm:locally finite lipschitz
    homology gives the locally finite homology}.]  
  We divide the proof into three steps:
  \begin{enumerate}
    \item For all~$L\in\bbR_{>0}$ and all~$K \in K(M)$, the
          inclusion~$C_*^{<L}(M, M - K) \rightarrow C_*(M, M
          - K)$ induces an isomorphism on homology.
    \item For all~$L \in \bbR_{>0}$, the inclusion~$\Clfl_*(M)
          \rightarrow \Clf_*(M)$ induces an isomorphism on homology.
    \item The inclusion~$\Clflip_*(M) \rightarrow \Clf_*(M)$ induces
          an isomorphism on homology.
  \end{enumerate}

  For the first step, let $L \in \bbR_{>0}$ and~$K \in K(M)$. We
  consider the commutative diagram
  \begin{align*}
    \xymatrix@=1.5em{
        0 
	\ar[r]
      & C_*^{<L}(M - K)
	\ar[r]\ar[d]
      & C_*^{<L}(M)
	\ar[r]\ar[d]
      & C_*^{<L}(M, M - K)
	\ar[r]\ar[d]
      & 0
        \\
	0
	\ar[r]
      & C_*(M - K)
	\ar[r]
      & C_*(M)
	\ar[r]
      & C_*(M, M - K)
	\ar[r]
      & 0
    }
  \end{align*}
  of chain complexes. By definition, the rows are exact; hence, there
  is a corresponding commutative diagram of long exact sequences in 
  homology. In view of the five
  lemma, it is therefore sufficient to show that the
  inclusion~$C_*^{<L}(U) \rightarrow C_*(U)$ induces an isomorphism on
  the level of homology whenever $U$ is an open subset of~$M$.

  By Lemma~\ref{lem:lipschitz homology the same} below, the
  inclusion~$\Clip_*(U) \rightarrow C_*(U)$ is a homology
  isomorphism. Let $\sd \colon C_*(U) \rightarrow C_*(U)$ be the
  barycentric subdivision operator. The map~$\sd$ is chain homotopic
  to the identity via a chain homotopy~$h \colon C_*(U) \rightarrow
  C_{*+1}(U)$~\cite{bredon}*{Section~IV.17}, and the classical
  construction of~$\sd$ and~$h$ shows that both~$\sd$ and~$h$ restrict
  to the Lipschitz chain complex~$\Clip_*(U)$. Moreover, for every
  Lipschitz simplex~$\sigma$ on~$U$ there is a~$k \in \bbN$ such
  that~$\lipschitz(\sd^k\sigma) < L$. Now the same argument as in the
  classical proof that singular homology is isomorphic to the homology
  of the chain complex of ``small''
  simplices~\cite{bredon}*{Section~IV.17} shows that the
  inclusion~$C_*^{<L}(U) \rightarrow \Clip_*(U)$ induces an
  isomorphism on homology. Therefore, $C_*^{<L}(U) \rightarrow C_*(U)$
  is a homology isomorphism.  This proves the first step.

  We now come to the proof of the second step. Since the structure maps in 
  the inverse system~$(C_*^{<L}(M, M - K))_{K \in
  K(M)}$ are surjective, we obtain a commutative diagram
  \begin{align*}
    \xymatrix@=1.5em{
        0 
	\ar[r]
      & \lim^1 H_{n+1}^{<L}(M,M - K)
        \ar[r]\ar[d]
      & \Hlfl_n(M)
        \ar[r]\ar[d]
      & \lim H_n^{<L}(M,M - K)
        \ar[d]\ar[r]
      & 0
        \\
        0
        \ar[r]
      & \lim^1 H_{n+1}(M,M - K)
        \ar[r]
      & \Hlf_n(M)
        \ar[r]
      & \lim H_n(M,M - K)
        \ar[r]
      & 0
      }
  \end{align*}
  with exact rows~\cite{weibel}*{Theorem~3.5.8}. By the first
  step, the outer vertical arrows are isomorphisms. Therefore, the
  five lemma shows that also the middle vertical arrows is an
  isomorphism, which proves the second step.

  Finally, the third step follows from the second step because
  homology is compatible with taking filtered colimits.
\end{proof}

\begin{lemma}\label{lem:lipschitz homology the same}
  Let $M$ be a Riemannian manifold and let $U \subset M$ be an open
  subset. Then the inclusion~$\Clip_*(U)
  \rightarrow C_*(U)$ induces an isomorphism on homology.
\end{lemma}
\begin{proof}
  The proof consists of an induction as, for example, in Bredon's
  proof of the de~Rham theorem~\cite{bredon}*{Section~V.9}:

  If $U \subset \bbR^n$ is a bounded convex subset, then one can easily
  construct a chain contraction for~$\Clip_*(U)$; therefore, the lemma
  holds for bounded convex subsets in Euclidean spaces. 

  If $U$,~$V \subset M$ are open subsets such that the lemma holds for
  both of them as well as for the intersection~$U \cap V$, then the
  lemma also holds for~$U \cup V$: The classical construction of
  barycentric subdivision (and the corresponding chain homotopy to the
  identity)~\cite{bredon}*{Section~IV.17} restricts to the Lipschitz
  chain complex and thus Lipschitz homology admits a Mayer-Vietoris
  sequence. 

  Proceeding by induction we see that the lemma holds for finite unions
  of bounded convex subsets of Euclidean space. Then a standard
  colimit argument shows that the lemma holds for arbitrary open
  subsets of Euclidean space.

  We call an open subset~$V$ of~$M$ \emph{admissible} if there is a
  smooth chart~$V' \rightarrow \bbR^n$ and a compact set~$K \subset M$
  such that $V \subset K \subset V'$. In particular, any admissible
  subset of~$M$ is bi-Lipschitz homeomorphic to an open subset
  of~$\bbR^n$, and hence the lemma holds for admissible subsets
  of~$M$. 

  Noting that the intersection of two admissible sets is admissible,
  the Mayer-Vietoris argument shows that the lemma holds for finite
  unions of admissible sets. Any open subset of~$M$ can be written as
  a union of admissible sets; hence, a standard colimit argument
  yields that the lemma holds for arbitrary open subsets of~$M$.
\end{proof}

\subsection{Cohomology with compact supports with a Lipschitz constraint}
\label{subsec:cslip cohomology}

The natural cohomological counterpart of locally finite homology is
cohomology with compact supports. Similarly, the cohomology theory
corresponding to Lipschitz locally finite homology is cohomology with
Lipschitz compact supports; here, ``corresponding'' means in
particular that there is an evaluation map linking homology and
cohomology (Remark~\ref{rem:lipschitz evaluation}).

\begin{definition}\label{def:lipschitz dual}
  Let $X$ be a metric space. A cochain~$f \in
  \hom_{\bbR}(\Clip_*(X),\bbR)$ is said to have \emph{Lipschitz
  compact support} if for all~$L \in \bbR_{>0}$ there exists a compact
  subset~$K \subset X$ such that
  \[ \fa{\sigma\in\map(\Delta^k,X)}
     \bigl( \lipschitz(\sigma) < L
            \land
	    \im(\sigma) \subset X - K
     \bigr)
     \Longrightarrow
     f(\sigma) = 0.
  \]
  The cochains with Lipschitz compact support form a subcomplex of the
  cochain complex~$\hom_{\bbR}(\Clip_*(X),\bbR)$; this subcomplex is
  denoted by~$\Ccslip^*(X)$.
  
  The cohomology of~$\Ccslip^*(X)$, denoted by~$\Hcslip^*(X)$, is
  called \emph{cohomology with Lipschitz compact supports}. 
\end{definition}

\begin{remark}\label{rem:lipschitz evaluation}
  Let $X$ be a metric space. By construction of the chain
  complexes~$\Clflip_*(X)$ and~$\Ccslip^*(X)$, the evaluation map
  \begin{align*}
    \langle\args,\args\rangle
    \colon 
    \Ccslip^*(X) \otimes \Clflip_*(X)
    & \longrightarrow
    \bbR \\
    f \otimes \sum_{i \in I} a_i \cdot \sigma_i
    & \longmapsto
    \sum_{i \in I} a_i \cdot f(\sigma_i)
  \end{align*}
  is well-defined. Moreover, the same computations as in the case of
  locally finite homology/cohomology with compact supports show that
  this evaluation descends to a map $\langle\args,\args\rangle \colon
  \Hcslip^*(X) \otimes \Hlflip_*(X) \longrightarrow \bbR$ on the level
  of (co)homology.
\end{remark}

Dually to Theorem~\ref{thm:locally finite lipschitz homology gives the
locally finite homology}, we obtain:

\begin{theorem}\label{thm:lipschitz dual same homology as compact support}
  For all connected Riemannian manifolds, the natural homomorphism~$\Ccs^*(M)
  \rightarrow \Ccslip^*(M)$ given by restriction induces an
  isomorphism on cohomology.
\end{theorem}
\begin{proof}
  We start by disassembling the cochain complex~$\Ccslip^*(M)$ into
  pieces that are accessible by the universal coefficient theorem:
  \begin{align*}
        \Ccslip^*(M) 
    & = \invlim_{L \rightarrow \infty} 
	\Ccsl^*(M), \\
	\Ccsl^*(M)
    & = \colim_{K \in K(M)}
	C^*_{<L}(M, M - K).
  \end{align*}
  Here, for all~$L \in \bbR_{>0}$ and all~$K \in K(M)$, 
  \begin{align*}  	 
          C^*_{<L}(M, M - K)
      &:= \hom_\bbR\bigl( C_*^{<L}(M, M - K), \bbR \bigr),\\
          \Ccsl^*(M)
      &:= \bigl\{ f \in C^*_{<L}(M, \emptyset)
                ;~i(f) \in \Ccs^*(M)
	  \bigr\},
  \end{align*}
  where $i \colon C^*_{<L}(M , \emptyset) \rightarrow C^*(M)$ is the
  map sending~$f$ to the extension~$\overline f$ of~$f$
  with~$\overline f(\sigma) = 0$ whenever $\lipschitz \sigma \geq L$. 

  Let $L \in \bbR_{>0}$ and $K \in K(M)$.  In the proof of
  Theorem~\ref{thm:locally finite lipschitz homology gives the locally
  finite homology}, we have shown that $C_*^{<L}(M,M - 
  K)\rightarrow C_*(M,M - K)$ induces an isomorphism on
  homology. Therefore, the restriction~$C^*(M, M - K)
  \rightarrow C_{<L}^*(M, M - K)$ induces an isomorphism on
  the level of cohomology by the universal coefficient theorem.
  Because homology commutes with colimits, it follows that the
  restriction map $\Ccs^*(M)\rightarrow \Ccsl^*(M)$ is a
  cohomology isomorphism.
 
  Notice that the structure maps in the 
  inverse system~$(\Ccsl^*(M))_{L \in \bbR_{>0}}$ are surjective, 
  in particular, they satisfy the Mittag-Leffler condition. 
  Furthermore, for~$L < L'$ there is a $k\in \bbN$ such
  that the $k$-fold barycentric subdivision~$\sd^k$ on~$C_*^{<L'}(M)$
  lands in~$C_\ast^{<L}(M)$. The classical construction of the
  barycentric subdivision operator shows that $\sd^k \colon
  C_*^{<L'}(M) \rightarrow C_*^{<L}(M)$ is a homotopy inverse of the
  inclusion~\cite{bredon}*{Section~IV.17}. Thus, the restriction map
  in cohomology $H^*_{\mathrm{cs},<L'}(M)\rightarrow \Hcsl^*(M)$ is
  surjective; in particular, the Mittag-Leffler condition on the
  level of cohomology is also satisfied. Therefore, the $\lim^1$-term
  vanishes~\cite{weibel}*{Proposition~3.5.7}, and we
  obtain~\cite{weibel}*{Theorem~3.5.8}
  \[       \Hcslip^*(M)
     \cong \invlim_{L \rightarrow \infty} \Hcsl^*(M)
     \cong \Hcs^*(M).
     \qedhere
  \]
\end{proof}

\subsection{Computing the Lipschitz simplicial volume via cohomology}
\label{subsec:lipvol via cohomology}

Any oriented, connected manifold possesses a (integral) 
\emph{fundamental class}, 
which is a distinguished generator of the locally finite homology 
$\Hlf_n(M;\bbZ)\cong\bbZ$ with integral coefficients 
in the top dimension $n=\dim(M)$. The fundamental class in 
$\Hlf_n(M)=\Hlf_n(M;\bbR)$ is, by definition, the image of the 
integral fundamental class under the coefficient change 
$\Hlf_n(M;\bbZ)\rightarrow\Hlf_n(M;\bbR)$. 
Correspondingly, 
one defines the \emph{cohomological} or \emph{dual fundamental class} 
as a distinguished generator of 
the top cohomology with compact supports.

\begin{definition}
  Let $M$ be an oriented, connected Riemannian \mbox{$n$-mani}\-fold
  (without boundary). The \emph{Lipschitz fundamental class} of~$M$ is the
  homology class~$[M]_\lip \in \Hlflip_n (M)$ that corresponds
  to the fundamental class~$[M] \in \Hlf_n(M)$ via the
  isomorphism~$\Hlflip_*(M) \rightarrow \Hlf_*(M)$
  (Theorem~\ref{thm:locally finite lipschitz homology gives
  the locally finite homology}). Analogously, one defines the 
  \emph{Lipschitz dual fundamental class}~$[M]^*_\lip \in \Hcslip^n(M)$ 
  of~$M$. 
\end{definition}

\begin{remark}\label{rem:integral comparisons}
The proofs of Theorem~\ref{thm:locally finite lipschitz homology gives
  the locally finite homology} and~\ref{def:lipschitz dual} work for
any coefficient module. Thus one can equivalently define the 
Lipschitz fundamental class as the image of the generator 
of $\Hlflip_n(M;\bbZ)$ that corresponds to the integral fundamental 
class in $\Hlf_n(M;\bbZ)\cong\Hlflip_n(M;\bbZ)$ under the change of 
coefficients $\bbZ\rightarrow\bbR$. Similar considerations apply to 
the Lipschitz dual fundamental class. 
\end{remark}

In the compact case, the simplicial volume can be expressed as the
inverse of the semi-norm of the dual fundamental
class~\cite{gromov}*{p.~17}. In the non-compact case, however, one has
to be a bit more
careful~\citelist{\cite{gromov}*{p.~17}\cite{loehphd}*{Theorem~C.2}}. Similarly,
also the Lipschitz simplicial volume can be computed in terms of
certain semi-norms on cohomology (Proposition~\ref{prop:simvol via
cohomology}).

\begin{definition}
  Let $M$ be a topological space, $k \in \bbN$, and let $A \subset
  \map(\Delta^k, M)$. 
  \begin{enumerate}
    \item For a locally finite chain~$c=\sum_{i \in I} a_i \cdot \sigma_i
          \in \Clf_k (M)$, let 
	  \[    \abs{c}_1^A
	     := \begin{cases}
                 \abs{c}_1 & \text{ if $\supp(c)\subset A$,}\\
                 \infty    & \text{ otherwise.}
                 \end{cases},
	  \]
          Here, $\supp(c) := \{i \in I ;~a_i \neq 0\}$.
    \item The semi-norms on (Lipschitz) locally finite/relative
          homology induced by~$\abs{\args}_1^A$ are denoted
          by~$\norm{\args}_1^A$.
    \item If $M$ is an oriented, connected $n$-manifold, then
          \[    \simvol M ^A 
             := \norm{[M]}^A_1. 
          \]
	  If moreover, $K \in K(M)$, then
	  \[    \simvol {M, M - K} ^A
	     := \norm{[M, M - K]}^A_1,
	  \]
	  where $[M, M -  K] \in H_n(M, M - K)$ is the
	  relative fundamental class. 
    \item For~$f \in C^k(M)$ we write
          \[    \supn f ^A
	     := \sup_{\sigma\in A} |f(\sigma)|
             \in [0, \infty].
	  \]
    \item The semi-norms on (relative) cohomology with (Lipschitz)
          compact supports induced by~$\supn{\cdot}^A$ are also
          denoted by~$\supn{\cdot}^A$. 
  \end{enumerate}
\end{definition}

\begin{proposition}[Duality principle for the Lipschitz simplicial volume]\label{prop:simvol via cohomology}
  Let $M$ be an oriented, connected Riemannian $n$-manifold. 
  \begin{enumerate}
    \item Then
          \[    \lipvol M
	     =  \inf\bigl\{ \simvol M ^A
	                  ;~A \in \liplfset_n (M)
		    \bigr\}.
	  \]
    \item Moreover, for all~$A \in \liplfset_n(M)$, we have
          \[    \simvol M ^A
	     =  \frac 1 {\supn{[M]^*_\lip}^A}.
	  \]
  \end{enumerate}
\end{proposition}
\begin{proof}
  The first part follows directly from the definitions. For the second
  part let $A \in \liplfset_n(M)$. Then 
  \begin{align*}
        \simvol M ^A
    & = \sup_{K \in K(M)} \simvol{M, M - K}^A \\
    & = \sup_{K \in K(M)} \frac 1 {\supn{[M, M - K]^*_\lip}^A}\\
    & = \frac 1 {\supn{[M]^*_\lip}^A}.
  \end{align*}
  We now explain these steps in more detail: 
  \begin{itemize}
    \item The first equality is shown by constructing an 
          appropriate diagonal sequence out of ``small'' relative
          fundamental cycles of the~$(M, M -  K)$ supported
          on~$A$~\cite{loehphd}*{Proposition~C.3}. 
    \item The class~$[M, M - K]^*_\lip \in
          H^n(\hom_\bbR(\Clip_*(M, M - K),\bbR))$ is the dual
          of the relative fundamental class in~$\Hlip_n(M, M - 
          K) \cong H_n(M, M - K) \cong \bbR$.

	  Therefore, the second equality is a consequence of the
          Hahn-Banach theorem -- this is exactly the same argument as
          in the non-Lipschitz case~\cite{loehphd}*{Proposition~C.6},
          but applied to functionals on~$\Clip_*(M)$ instead
          of~$C_*(M)$; this is possible because $A$ is Lipschitz.
    \item The last equality is equivalent to 
          \[\inf_{K\in K(M)}\supn{[M, M - 
            K]^*_\lip}^A=\supn{[M]^*_\lip}^A.
          \]
          Here the $\ge$-inequality is clear. For the
          $\le$-inequality, let $\epsilon>0$ and
          consider~$f\in\Ccslip^n(M)$ with $\supn
          f^A\le\supn{[M]_\lip^\ast}^A+\epsilon$. By
          Theorem~\ref{thm:lipschitz dual same homology as compact
            support}, there is a compactly supported cochain $g$ and a
          $(n-1)$-cochain~$h$ with Lipschitz compact support such that
          $f=g+\delta h$. Since $A\in\liplfset_n (M)$, the chain~$h'$
          defined by
          \[h'(\sigma)=\begin{cases}
                      h(\tau) &\text{ if $\sigma\in
                        \bigcup_{j=0}^n\{\partial_j\sigma;~\sigma \in A~\}$, }\\ 
                      0         &\text{ otherwise,} 
                      \end{cases}
          \]
          is compactly supported. Further, 
          $f'\defq g+\delta h'$ is compactly supported, cohomologous
          in~$\Ccslip^*(M)$ to~$f$,
          and $\supn{f'}^A=\supn{f}^A$. In particular, there is $K\in
          K(M)$ with $f'\in C^n_\lip(M,M - K)$ and 
          \[\supn{ [M,M - K]_\lip^\ast}^A\le\supn{ f'}^A=\supn{
            f}^A\le\supn{[M]_\lip^\ast}^A+\epsilon.\]
 \end{itemize}
  This finishes the proof of the duality principle.
\end{proof}

\subsection{Product structures in the Lipschitz setting}
\label{subsec:cross-products}

The definition of product structures in singular (co-)homology is
based on the following maps: Let $X$ and~$Y$ be topological
spaces. Then there exist chain maps $\EZ: C_\ast(X)\otimes
C_\ast(Y)\rightarrow C_\ast(X\times Y)$ and $\AW: C_\ast(X\times
Y)\rightarrow C_\ast(X)\otimes C_\ast(Y)$, called the
\emph{Eilenberg-Zilber map} and the \emph{Alexander-Whitney map},
respectively, such that $\EZ\circ \AW$ and $\AW\circ \EZ$ both are
naturally homotopic to the identity; explicit formulas are, for
example, given in Dold's book~\cite{dold}*{12.26 on p.~184}.  

The map $\EZ$ and the composition~$C^\ast(X)\otimes
C^\ast(Y)\rightarrow C^\ast(X\times Y)$, $f\otimes g\mapsto (f\otimes
g)\circ \AW$, induce the so-called \emph{cross-products}
\begin{gather}
\times\colon H_m(X)\otimes H_n(Y)\rightarrow H_{m+n}(X\times Y)\label{eq:cross homology},\\
\times\colon H^m(X)\otimes H^n(Y)\rightarrow H^{m+n}(X\times Y)\nonumber
\end{gather}
in homology and cohomology, respectively. 

Next we describe these cross-products more explicitly on the (co)chain
level: Let $f\in C^m(X)$ and $g\in C^n(Y)$. 
Let $\pi_X$ and $\pi_Y$ be the projections from $X\times Y$ to 
$X$ and $Y$, respectively. For a $k$-simplex $\sigma$, 
let $\frface l \sigma$ and $\baface {(k-l)} \sigma$ the 
\emph{$l$-front face} and the \emph{$(k-l)$-back face} of~$\sigma$, 
respectively.  
Then the explicit formula for~$\AW$ in \emph{loc.~cit.} yields 
\begin{equation}\label{eq:explicit cohomological cross product}
(f\times g)(\sigma)=f(\pi_X\circ\frface m \sigma )\cdot g(\pi_Y\circ\baface n \sigma ). 
\end{equation}
For simplices $\sigma:\Delta^m\rightarrow X$ and
$\rho:\Delta^n\rightarrow Y$, the chain $\EZ(\sigma\otimes\rho)$ can
be described as follows: The
product~$\Delta^n\times\Delta^m\rightarrow X\times Y$ of $\sigma$ and
$\rho$ is not a simplex but can be chopped into a union of
$(m+n)$-simplices (like a square can be chopped into triangles, or a
prism into tetrahedra). Then $\EZ(\sigma\otimes\rho)$ is the sum of
these $(m+n)$-simplices.

From this description we see that 
if $c=\sum_i a_i\sigma_i$ and $d=\sum_j b_j\rho_j$ 
are (Lipschitz) locally finite chains in (metric) spaces 
$X$ and $Y$, 
then $\sum_{i,j}a_ib_j(\sigma_i\times\rho_j)$ is a 
(Lipschitz) locally finite chain in $X\times Y$. 
Thus,~\ref{eq:cross homology} extends to maps 
\begin{gather*}
\times\colon\Hlf_m(X)\otimes\Hlf_n(Y)\rightarrow\Hlf_{m+n}(X\times Y),\\
\times\colon\Hlflip_m(X)\otimes\Hlflip_n(Y)\rightarrow\Hlflip_{m+n}(X\times Y).
\end{gather*}

In general, the cross-product of two cocycles with compact supports
has not necessarily compact support. However, the cross-product of two
cochains with Lipschitz compact supports again has Lipschitz compact
support:

\begin{lemma}\label{lem:cross product restricts}
  Let $M$ and $N$ be two complete metric spaces, and let $m$,~$n \in
  \bbN$. Then the cross-product on $C^*(M) \otimes C^*(N) \rightarrow
  C^*(M\times N)$ restricts to a cross-product
  \[ \times \colon 
     \Ccslip^m(M)\otimes\Ccslip^n(N)
     \rightarrow
     \Ccslip^{m+n}(M \times N), 
  \]
  which induces a cross-product~$\Hcslip^m(M)\otimes\Hcslip^n(N)
  \rightarrow \Hcslip^{m+n}(M \times N)$. 
\end{lemma}
\begin{proof}
  Let $f\in\Ccslip^m(M)$ and $g\in\Ccslip^n(N)$. 
  Let $L \in \bbR_{>0}$. Because $f$ and $g$ are cochains with
  Lipschitz compact supports, there are compact sets~$K_M \subset M$
  and $K_N \subset N$ with 
  \[ \fa{\sigma \in \map(\Delta^m, M)}
     \bigl(\lipschitz(\sigma) \leq L
           \;\land\;
           \im (\sigma) \subset M - K_M
     \bigr)
     \Longrightarrow f(\sigma) = 0,
  \]
  and analogously for~$g$ and~$K_N$.

  We now consider the compact set~$K := U_L(K_M) \times U_L(K_N)
  \subset M\times N$, where $U_L(X)$ denotes the set of all points
  with distance at most~$L$ from~$X$. Because the diameter of the
  image of a Lipschitz map on a standard simplex is at most as large
  as $\sqrt 2$ times the Lipschitz constant of the map in question, we
  obtain: If $\sigma \in \map(\Delta^{m+n}, M \times N)$ with
  $\lipschitz(\sigma) \leq L$ and $\im(\sigma) \subset M \times N
  - K$, then 
  \[ \im(\pi_M \circ \sigma) \subset M - K_M
     \qquad\text{or}\qquad
     \im(\pi_N \circ \sigma) \subset N - K_N.
  \]
  In particular, $f(\pi_M \circ \frface m \sigma) = 0$ or $g (\pi_N \circ
  \baface n \sigma) = 0$. By~\ref{eq:explicit cohomological cross product}, 
  $(f \times g) (\sigma) =0$. In
  other words, the cross-product~$f \times g$ lies in~$\Ccslip^{m+n}(M
  \times N)$.
\end{proof}

\begin{definition}\label{def:projection sets of simplices}
  Let $M$ and $N$ be two topological spaces, let $m$,~$n\in\bbN$, and
  let $A \subset \map(\Delta^{m+n}, M\times N)$. Then we write
  \begin{align*}
        A_M 
   & := \bigl\{ \pi_M \circ \frface m \sigma
	      ;~\sigma \in A
	\bigr\},
        \\ 
        A_N 
   & := \bigl\{ \pi_N \circ \baface n \sigma
	      ;~\sigma \in A
	\bigr\},
  \end{align*}
  where $\pi_M \colon M\times N \rightarrow M$ and $\pi_N \colon M \times
  N \rightarrow N$ are the projections. Notice that $A_M$ and~$A_N$
  depend on~$m$ and~$n$, but the context will always make clear which
  indices are involved; therefore, we suppress~$m$ and~$n$ in the
  notation. 
\end{definition}

The cross-product of cochains with Lipschitz compact support is
continuous in the following sense:

\begin{remark}\label{rem:products and norms}
  Let $M$ and $N$ be two topological spaces, and let~$m$,~$n \in
  \bbN$. Then -- by the explicit 
  description~\ref{eq:explicit cohomological cross product} -- 
  the cross-product satisfies 
  \[      \supn {f \times g}^A
     \leq \supn f ^{A_M}
          \cdot \supn g ^{A_N}
  \]
  for all~$A \subset \map(\Delta^{m+n}, M\times N)$ and all~$f \in
  \Ccslip^m (M)$, $g \in \Ccslip^n (N)$.
\end{remark}

Notice however, that in general the sets~$A_M$ and~$A_N$ are
\emph{not} locally finite even if $A$ is locally finite. This issue is
addressed in Section~\ref{subsec:sparse}.

\begin{lemma}\label{lem:product fundamental cocycles}
  Let $M$ and $N$ be oriented, connected, complete Riemannian
  manifolds. Then
  \[ [M \times N]^*_\lip = [M]^*_\lip \times [N]^*_\lip
     \in \Hcslip^{*}(M \times N).
  \]
\end{lemma}
\begin{proof}
In view of Remark~\ref{rem:integral comparisons}, it is 
enough to show 
\[\bigl\langle [M]^*_\lip \times [N]^*_\lip, [M]_\lip\times
[N]_\lip\bigr\rangle=1.\] Let $f\in\Ccslip^m(M)$ and $g\in\Ccslip^n(N)$ be
fundamental cocycles that vanish on degenerate simplices. Such
fundamental cocycles always exist; for example, let $f$ be the cocycle
$\sigma\mapsto\int_{\Delta^m}\sigma^\ast\omega$ where
$\omega\in\Omega^m(M)$ is a compactly supported differential $m$-form
with $\int_M\omega=1$. Note that the integral exists by Rademacher's
theorem~\cite{evans+gariepy}.  Let $w=\sum_i
a_i\sigma_i\in\Clflip_m(M)$ and $z=\sum_j b_j\rho_j\in\Clflip_n(N)$ be
fundamental cycles of~$M$ and $N$ respectively.  The Eilenberg-Zilber
and Alexander-Whitney maps have the property that $\AW \circ \EZ$
differs from the identity by degenerate
chains~\cite{eilenberg}*{Theorem~2.1a (2.3)}.  Thus, we obtain
\begin{align*}
\langle f\times g, w\times z\rangle 
 &=\sum_{i,j}a_ib_j (f\times g)(\sigma_i\times \rho_j)\\
 &=\sum_{i,j}a_ib_j(f\otimes g)\bigl(\AW\circ \EZ (\sigma_i\otimes \rho_j)\bigr)\\
 &=\sum_{i,j}a_ib_j(f\otimes g)\bigl(\sigma_i\otimes\rho_j+\text{degenerate~simplices}\bigr)\\
 &=\sum_{i,j}a_ib_jf(\sigma_i)g(\rho_j)=f(w)g(z)=1.\qedhere
\end{align*}
\end{proof}

\subsection{Representing the fundamental class of the product by sparse cycles}\label{subsec:sparse}

The functor~$\Clf_*$ is only functorial with respect to \emph{proper}
maps. For example, in general, the projection of a locally finite
chain on a product of non-compact spaces to one of its factors is not
locally finite.

\begin{definition}
  Let $M$ and $N$ be two topological spaces, and let $k\in \bbN$. A
  locally finite set~$A \in \lfset_k(M\times N)$ is called
  \emph{sparse} if 
  \[ \bigl\{ \pi_M \circ \sigma
           ;~\sigma \in A
     \bigr\}
     \in \lfset_k(M)
     \qquad\text{and}\qquad
     \bigl\{ \pi_N \circ \sigma
           ;~\sigma \in A
     \bigr\}
     \in \lfset_k(N), 
  \]
  where $\pi_M \colon M\times N \rightarrow M$ and $\pi_N \colon M \times
  N \rightarrow N$ are the projections.
  
  A locally finite chain~$c \in \Clf_*(M \times N)$ is called
  \emph{sparse} if its support is sparse.
\end{definition}

The following proposition is crucial in proving the product inequality 
for the Lipschitz simplicial volume. 

\begin{proposition}\label{prop:sparse cycles}
  Let $M$ and $N$ be two oriented, connected, complete Riemannian
  manifolds (without boundary) with non-positive sectional curvature.
  \begin{enumerate}
    \item
      For any cycle~$c \in \Clflip_*(M \times N)$ there is a
      sparse cycle~$c' \in \Clflip_*(M\times N)$ satisfying
      \[ \abs{{c'}}_1 \leq \abs{c}_1
         \qquad\text{and}\qquad
         c \sim c' \text{\ in~$\Clflip_*(M \times N)$}.
      \]
    \item
      In particular, the Lipschitz simplicial volume can be computed
      via sparse fundamental cycles, i.e., 
      \[ \lipvol{M \times N}
         = \inf 
	   \bigl\{ \norm{M\times N}^A
	         ;~A \in \liplfset_{\dim M + \dim N}(M \times N),\ 
	           \text{$A$ sparse}
	   \bigr\}
         .
      \]
  \end{enumerate}
\end{proposition}

\begin{proof}
  The second part is a direct consequence of the first part. For the
  first part, we take advantage of a straightening procedure:

  Let $F_M \subset M$ and $F_N \subset N$ be locally finite subsets
  with $U_1(F_M) = M$ and $U_1(F_N) = N$. Then the corresponding
  preimages~$\ucov F_M := p_M^{-1}(F_M) \subset \ucov M$ and~$\ucov
  F_N := p_N^{-1}(F_N) \subset \ucov N$ satisfy~$U_1(\ucov F_M)
  = \ucov M$ and~$U_1(\ucov F_N) = \ucov N$, where $p_M \colon \ucov
  M \rightarrow M$ and $p_N \colon \ucov N \rightarrow N$ are the
  Riemannian universal covering maps.

  Furthermore, also the product~$F := F_M \times F_N \subset M \times
  N$ is locally finite, and there is a $\pi_1(M) \times
  \pi_1(N)$-equivariant map~$f \colon \ucov M \times \ucov N
  \rightarrow \ucov F_M \times \ucov F_N =: \ucov F$ such that
  \[ d_{\ucov M \times \ucov N}\bigl(z, f(z)\bigr) 
    \leq \sqrt 2
  \]
  holds for all~$z \in \ucov M \times \ucov N$. 

  For~$\sigma \in \map(\Delta^k,M\times N)$, we define
  \[      h_\sigma 
      :=  (p_{M} \times p_N) \circ 
          \bigl[ \ucov \sigma, 
                 [f(\ucov \sigma(v_0)), \dots, f(\ucov \sigma(v_k))]
          \bigr]
          \colon \Delta^k \times [0,1] 
          \rightarrow M \times N, 
  \]
  where $v_0, \dots, v_k$ are the vertices of the standard
  simplex~$\Delta^k$, and $\ucov \sigma$ is a lift of $\sigma$. 
  By Lemma~\ref{lem:cat0 estimate} and
  Proposition~\ref{prop:geodesic simplex - non-positive curvature}
  (and Remark~\ref{rem:diameter of geodesic simplices}),
  the map~$h_\sigma$ is Lipschitz, and the Lipschitz constant can be
  estimated from above in terms of the Lipschitz constant
  of~$\sigma$. 
  Moreover, the fact that $f$ is equivariant and covering theory show
  that 
  \begin{align}\label{eq:h is compatible with boundaries}
      h_{\sigma \circ \partial_j}
    = h_\sigma \circ (\partial_j \times \id_{[0,1]})
  \end{align}
  for all~$\sigma \in \map(\Delta^k,M \times N)$ and all~$j \in \{0,
  \dots, k-1\}$, where $\partial_j \colon \Delta^{k-1} \rightarrow
  \Delta^k$ is the inclusion of the $j$-th face.

  We now consider the map
  \begin{align*}
      H \colon \Clflip_*(M \times N)
    & \longrightarrow
      \Clflip_{*+1}(M \times N) \\
      \sum_{i \in I} a_i \cdot \sigma_i
    & \longmapsto
      \sum_{i \in I} a_i \cdot \overline h_{\sigma_i},
  \end{align*}
  where $\overline h_{\sigma}$ is the singular chain constructed out
  of~$h_{\sigma}$ by subdividing the prism~$\Delta^k\times [0,1]$ in
  the canonical way into a sum of $k+1$ simplices of dimension~$k+1$
  (compare Lemma~\ref{lem:lemma from lee}). 
  
  The map~$H$ is indeed well-defined: As discussed above, 
  for all~$c \in \Clflip_k(M \times N)$, all simplices occurring
  in the (formal) sum~$H(c)$ satisfy a uniform Lipschitz condition
  depending on~$\lipschitz(c)$. Further, it follows from 
  $\im(h_\sigma)\subset U_{\sqrt{2}}(\im(\sigma))$ that 
  $H$ maps locally finite chains to locally finite chains. 
  As next step, we define 
  \begin{align*}
      \phi \colon \Clflip_*(M \times N)
    & \longrightarrow 
      \Clflip_*(M \times N)
      \\
      \sum_{i \in I} a_i \cdot \sigma_i
    & \longmapsto
      \sum_{i \in I} a_i \cdot h_{\sigma_i}(\args, 1).
  \end{align*}
  In other words, $\phi$ is given by replacing each simplex by a
  straight simplex whose vertices lie in~$F_M \times F_N$ and whose
  vertices are close to the ones of the original
  simplex. Property~\ref{eq:h is compatible with boundaries} implies 
  that $\varphi$ is a chain map and that $H$ is a
  chain homotopy between the identity and~$\varphi$ (see 
  Lemma~\ref{lem:lemma from lee}).

  By construction, $\norm \phi \leq 1$. Therefore, it remains to show
  that the image of~$\varphi$ contains only sparse chains: 

  Let $c \in \Clflip_k(M \times N)$. Let~$A := \supp(\varphi(c))$.
  Because the geodesics in~$\ucov M \times \ucov N$ are just products
  of geodesics in~$\smash{\ucov M}$ and~$\smash{\ucov N}$, it follows
  that the projection~$\pi_M \colon M\times N \rightarrow M$ preserves
  straight simplices.  Thus, the set $\bigl\{ \pi_M \circ
  \sigma;~\sigma \in A\bigr\}$ consists of straight simplices whose
  Lipschitz constant is bounded by~$\lip(c)$ and whose vertices lie
  in~$F_M$.  The fact that $F_M$ is locally finite and that there are
  only finitely many straight simplices with a bounded Lipschitz
  constant and the same vertices imply that~$\{\pi_M \circ
  \sigma;~\sigma \in A\}$ is locally finite.  Similarly for the
  projection to~$N$. So the chain~$\varphi(c)$ is sparse.
\end{proof}

\subsection{Conclusion of the proof of the product inequality}
\label{subsec:product formula conclusion}

Finally, we can put all the pieces collected in the previous sections
together to give a proof of the product inequality:

\begin{proof}[Proof of Theorem~\ref{thm:product formula}.]
   In the following, we write $m := \dim M$ and $n := \dim N$. In
   order to prove the product inequality, it suffices to find for
   each~$\epsilon \in \bbR_{>0}$ locally
   finite sets~$A_M \in \liplfset_m(M)$ and~$A_N \in \liplfset_n(N)$
   with
   \[      \simvol M ^{A_M} \cdot \simvol N ^{A_N}
      \leq \lipvol{M\times N} + \epsilon.
   \]
   
   For every~$\epsilon\in\bbR_{>0}$, Proposition~\ref{prop:sparse cycles}
   provides us with a sparse fundamental cycle~$c \in \Clflip_{m+n}(M
   \times N)$ with support $A$ satisfying
   $\abs c _1 \leq \lipvol{M \times N} + \epsilon$.
   In particular, 
   \[ \lipvol{M\times N}^A\leq \lipvol{M \times N} + \epsilon. 
   \]
   By sparseness, the sets 
   $A_M$ and $A_N$ associated to 
   $A$ (see Definition~\ref{def:projection sets of simplices}) lie in 
   $\liplfset_m(M)$ and $\liplfset_n(N)$, respectively. 
   The duality principle (Proposition~\ref{prop:simvol via
   cohomology}) yields
   \[   \simvol M ^{A_M}
      = \frac 1 {\supn{[M]^*_\lip}^{A_M}},
      \quad
        \simvol N^{A_N}
      = \frac 1 {\supn{[N]^*_\lip}^{A_N}},
      \quad
        \simvol {M \times N}^A
      = \frac 1 {\supn{[M \times N]^*_\lip}^A};
   \]
   the cohomological terms are related as follows
   \[       \supn{[M \times N]^*_\lip}^A 
       \leq       \supn{[M]^*_\lip}^{A_M}
            \cdot \supn{[N]^*_\lip}^{A_N}
   \]
   because~$[M \times N]^*_\lip = [M]^*_\lip \times [N]^*_\lip$
   (Lemma~\ref{lem:product fundamental cocycles}) and the
   cohomological cross-product is compatible with the semi-norms
   (Remark~\ref{rem:products and norms}).

   Therefore, we obtain
   \[
             \lipvol{ M} \cdot \lipvol{N}
      \leq  \simvol M^{A_M} \cdot \simvol N^{A_N}\\
      \leq  \simvol {M \times N}^A\\
      \leq  \lipvol{M \times N} + \epsilon.\qedhere
   \]
\end{proof}

\section{Proportionality principle for non-compact
  manifolds}\label{sec:proportionality for non-compact manifolds} 

\noindent 
Thurston's proof of the proportionality principle in the compact case
is based on ``smearing'' singular chains to so-called measure
chains~\citelist{\cite{thurston}*{p.~6.6--6.10}\cite{loehdiplom}*{Chapter~5}}.
We prove the proportionality principle in the non-compact
case by combining the smearing technique with a discrete approximation 
of it; to this end, we replace measure homology by Lipschitz
measure homology, a variant that incorporates a Lipschitz constraint
(Section~\ref{subsec:smearing map}).

Throughout Section~\ref{sec:proportionality for non-compact
manifolds}, we often refer to the following setup:

\begin{setup}\label{setup:smearing map}
Let $M$ and $N$ be oriented, connected, complete, non-positively curved
Riemannian manifolds of finite volume without boundary whose
universal covers are isometric. We denote the common universal
cover by~$U$. Let $G=\isom^+(U)$ be its group of orientation-preserving
isometries. Then $\Gamma=\pi_1(M)$ and $\Lambda=\pi_1(N)$ are lattices
in $G$ by Lemma~\ref{lem:fundamental group is lattice} below.  Let
$\mu_{\Lambda\bs G}$ denote the normalized Haar measure on $\Lambda\bs
G$. The universal covering maps of $M$ and $N$ are denoted by $p_M$
and $p_N$, respectively.
\end{setup}

The following lemma is well known for locally symmetric spaces 
and compact manifolds 
but we were unable to find a reference in the general case. 

\begin{lemma}\label{lem:fundamental group is lattice}
Let $M$ be a complete Riemannian manifold of finite volume. Then 
$\Gamma=\pi_1(M)$ is a lattice in $G=\isom(\ucov M)$. 
\end{lemma}

\begin{proof}
The isometry group $G$ acts smoothly and properly on 
$\ucov M$. It is easy to see that $\Gamma$ is a discrete 
subgroup. 
Let $x_0\in\ucov M$, and let $K\subset G$ be the stabilizer 
of $x_0$. 
Let $\nu\rightarrow Gx_0$ be the normal bundle of $Gx_0$, and let 
$\nu(r)$ denote the sub-bundle of vectors of length at most~$r$. 
By the slice 
theorem~\citelist{\cite{palais}*{Section~2.2}\cite{duistermaat}*{Chapter~2}}, 
there exists $r>0$ such 
that the exponential map $\exp:\nu(r)\rightarrow V$ is a diffeomorphism 
onto a tubular neighborhood $V$ of $Gx_0$. The map 
$f:G\times_K \nu_{x_0}\rightarrow\nu$, $(g,z)\mapsto Tg(z)$ is a 
diffeomorphism. Define $g=f^{-1}\circ\exp^{-1}$. We equip 
$G/K$ with the Riemannian metric that turns the diffeomorphism 
$G/K\rightarrow Gx_0$ into an isometry. Since $\nu_{x_0}$ can be 
equipped with a  
$K$-invariant metric ($K$ is compact), it is easy to see that 
$G\times_K\nu_{x_0}(r)$ carries a $G$-invariant Riemannian metric such 
that the projection $G\times_K\nu_{x_0}(r)\rightarrow G/K$ is a 
Riemannian submersion. 
By compactness, there is $\lambda>0$ such that 
$T_zg$ has norm at most $\lambda$ for all $z\in\exp(\nu_{x_0}(r))$. 
By $G$-invariance of the metrics, $T_z g$ has norm at most $\lambda$ 
for all $z\in V$, thus, $g$ is $\lambda$-Lipschitz, and so is 
the induced map between the $\Gamma$-quotients. 
We obtain that 
\[\vol\bigl(\Gamma\bs (G\times_K\nu_{x_0}(r))\bigr)\le
\lambda^{\dim(M)}\vol(\Gamma\bs V)<\infty.\] 
Fubini's theorem for Riemannian 
submersions~\cite{sakai}*{Theorem~5.6 on p.~66} yields 
\[\vol(\Gamma\bs G/K)\vol(\nu_{x_0}(r))=
\vol\bigl(\Gamma\bs (G\times_K\nu_{x_0}(r))\bigr)<\infty.\]
Thus $\vol(\Gamma\bs G/K)<\infty$. Now equip $G$ with a $G$-equivariant 
metric such that $G\rightarrow G/K$ is a Riemannian submersion. 
By uniqueness, the corresponding Riemannian measure on $G$ is a 
Haar measure. Fubini's theorem and $\vol(\Gamma\bs G/K)<\infty$ 
show that $\vol(\Gamma\bs G)<\infty$. 
\end{proof}

\subsection{Integrating Lipschitz chains}

Before introducing the smearing operation in
Section~\ref{subsec:smearing map}, we first discuss integration of
Lipschitz chains, which provides a means to detect which class in
locally finite homology a given Lipschitz cycle represents.

Let $M$ be an $n$-dimensional Riemannian manifold, and let $K\subset
M$ be a compact, connected subset with non-empty interior. 
Let $\Omega^\ast(M,M-K)$ be the kernel of the restriction
homomorphism 
$\Omega^\ast(M)\rightarrow\Omega^\ast(M-K)$ on differential
forms. The corresponding cohomology groups are denoted by 
$\Hde^\ast(M,M-K)$. The de~Rham map $\Omega^\ast(M)\rightarrow
C^\ast(M)$ restricts to the respective kernels and thus induces a
homomorphism, called \emph{relative de Rham map}, 
\[\Psi^\ast:\Hde^\ast(M,M-K)\rightarrow H^\ast(M,M-K).\]
The relative de Rham map is an isomorphism, which follows from the
bijectivity of the absolute de Rham map and an application of the
five~lemma. Note that integration gives a homomorphism
$\int:\Hde^n(M,M-K)\rightarrow\bbR$. Moreover, it is well known that 
\begin{equation}\label{eq:compare relative evaluation with integration}
\bigl\langle \Psi^n[\omega], [M,M-K]\bigr\rangle=\int_M\omega
\end{equation}
holds for all $n$-forms~$\omega$.

\begin{proposition}\label{prop:evaluation at the fundamental class}
  Let $M$ be a Riemannian $n$-manifold, and let $c=\sum_{k\in
  \bbN}a_k\sigma_k \in \Clf_n(M)$ be a cycle with~$\abs c _1 < \infty$
  and $\lipschitz(c) < \infty$.
  \begin{enumerate}
    \item Then $\langle\dvol_M,\sigma_k\rangle \leq \lipschitz(c)^n\vol(\Delta^n)$
          for every~$k\in\bbN$.
    \item Furthermore, we have the following equivalence:
          \begin{equation*}
            \sum_{k\in\bbN} a_k \cdot \langle\dvol_M,\sigma_k\rangle = \vol(M)
            \iff
            \text{$c$ is a fundamental cycle.}
          \end{equation*}
  \end{enumerate}
\end{proposition}

\begin{proof}
  For the first part, it suffices to observe that all Lipschitz
  simplices~$\sigma$ are almost everywhere differentiable, that
  $\sigma^* \dvol_M$ is measurable (by Rademacher's
  theorem~\cite{evans+gariepy}), and that  
  \begin{equation*}
       \absfix{\langle \dvol_M,\sigma\rangle}
    =  \biggl|\int_{\Delta^n}\sigma^\ast \dvol_M \biggr|
   \le \esssup_{x\in\Delta^n}\norm{T_x\sigma}^n\vol\bigl(\Delta^n\bigr)
   \le \lipschitz(\sigma)^n\vol\bigl(\Delta^n\bigr)
  \end{equation*}
  holds. In particular, we see that $\sum_{k\in\bbN} a_k\langle
  \dvol_M, \sigma_k\rangle$ converges absolutely.

  For the second part, let $s \in \bbR$ be the number defined by 
  \[ [c] = s \cdot [M] \in \Hlf_n(M).
  \]
  In the following, we show that $\sum_{k\in\bbN} a_k\langle\dvol_M,
  \sigma_k\rangle = s \cdot \vol(M)$: To this end, we first relate~$s
  \cdot \vol(K)$ for compact~$K$ to a finite sum derived from the
  series on the left hand side, and then use a limit process to
  compute the value of the whole series. 

  Let $K \subset M$ be a connected, compact subset with non-empty
  interior. For~$\delta \in \bbR_{>0}$ let $g_\delta \colon M
  \rightarrow [0,1]$ be a smooth function supported on the closed
  \mbox{$\delta$-neigh}\-bor\-hood~$K_\delta$ of~$K$ with~$g|_K =
  1$. Then $g_\delta \cdot \dvol_M \in \Omega^n(M, M - K_\delta)$ is a
  cocycle, and
  \begin{align*}
      s \cdot \vol(K) 
    = \lim_{\delta \rightarrow 0} s \cdot \int_M g_\delta \dvol_M. 
  \end{align*}
  On the other hand, the map~$H_n(j_\delta) \colon \Hlf_n(M) \rightarrow
  H_n(M, M - K_\delta)$ induced by the inclusion~$j_\delta \colon (M,
  \emptyset) \rightarrow (M, M-K_\delta)$ maps the fundamental class
  of~$M$ to the relative fundamental class of~$(M, M-K_\delta)$ and
  $H_n(j_\delta)[c]$ is represented by~$\sum_{\im \sigma_k \cap K_\delta \neq
  \emptyset} a_k \sigma_k$.  Therefore, we obtain 
  by~\ref{eq:compare relative evaluation with integration}
  \begin{align*}    
        \lim_{\delta\rightarrow 0}
        \sum_{\im \sigma_k \cap K_\delta \neq \emptyset}
        a_k \cdot \langle g_\delta \cdot \dvol_M, \sigma_k \rangle
    &=  \lim_{\delta\rightarrow 0}
        \bigl\langle
          \Psi^n[g_\delta \cdot \dvol_M], s \cdot[M, M-K_\delta]
        \bigr\rangle\\
    &=  \lim_{\delta\rightarrow 0}
        s \cdot \int_M g_\delta \dvol_M\\ 
    &=  s \cdot \vol(K).
  \end{align*}
  For each~$k \in \bbN$ and~$\delta \in
  \bbR_{>0}$ we have
  $  | \langle g_\delta \cdot \dvol_M, \sigma_k
            \rangle
     |
     \leq \lipschitz(c)^n \vol(\Delta^n)$, 
  and hence
  \begin{multline*}
    \biggl| \sum_{k \in\bbN}
            a_k \langle \dvol_M, \sigma_k 
                \rangle
          - \sum_{\im\sigma_k \cap K_\delta\neq 0}
            a_k \langle g_\delta\cdot\dvol_M, \sigma_k
                \rangle
    \biggr| \\
    \leq 2 \lipschitz(c)^n \vol(\Delta^n)
         \cdot \sum_{\im\sigma_k \subset M- K}
               |a_k|.
  \end{multline*}
  
  Because $\sum_{k \in \bbN} |a_k| < \infty$, there is an exhausting
  sequence~$(K^m)_{m \in \bbN}$ of compact, connected subsets of~$M$
  with non-empty interior satisfying
  \[ \lim_{m \rightarrow \infty} \vol(K^m) = \vol(M) 
     \qquad
     \text{and}
     \qquad
     \lim_{m \rightarrow \infty} 
     \sum_{\im \sigma_k \subset M - K^m}
          |a_k|
     = 0.
  \]
  Thus, the estimates of the previous paragraphs yield
  \begin{align*}
        \sum_{k\in\bbN} a_k \cdot 
                            \langle \dvol_M,\sigma_k \rangle
     &= \lim_{m \rightarrow \infty}
	\lim_{\delta \rightarrow 0}
        \sum_{\im\sigma_k \cap K^m_\delta \neq \emptyset} 
	     a_k \cdot \langle g_\delta^m\cdot \dvol_M, \sigma_k
                       \rangle \\
     &= \lim_{m \rightarrow \infty}
        s \cdot \vol(K^m)\\
     &= s \cdot \vol(M).  
  \end{align*}

  If $c$ is a fundamental cycle, then $s=1$ and hence the series has
  value~$\vol(M)$. Conversely, if the series evaluates to~$\vol(M)$,
  then $\vol(M)$ must be finite by the first part. Therefore, we can
  deduce from the computation above that~$s=1$, i.e., $c$ is a
  fundamental cycle. 
\end{proof}

One should be aware that the (locally finite) simplicial volume of a
non-compact manifold~$M$ might be finite even if $\vol(M) = \infty$, 
e.g., $\lfvol{\bbR^2}=0$ -- unlike the Lipschitz simplicial volume  
as the following direct corollary of 
Proposition~\ref{prop:evaluation at the fundamental class} shows.

\begin{corollary}\label{cor:finite vol and finite lipvol}
  Let $M$ be a Riemannian manifold. If $\lipvol M$ is finite, then so
  is~$\vol(M)$.
\end{corollary}

\subsection{The smearing homomorphism}\label{subsec:smearing map}

Let $M$ and $N$ be smooth manifolds (with or without boundary). 
The set of smooth maps $M\rightarrow N$ equipped with the topology 
that turns the differential map from this set to $\map(TM,TN)$ into 
a homeomorphism onto its image is denoted by $C^1(M,N)$. This topology
is called \emph{$C^1$-topology}.  

The following defines a variant of Thurston's measure 
homology~\cite{thurston}*{p.~6.6f}. 
\begin{definition}[Lipschitz measure homology]\label{def:variants of
    measure homology} 
Let $M$ be a Riemannian manifold. 
\begin{enumerate}\renewcommand{\theenumi}{\alph{enumi}}
\item 
A signed Borel measure $\mu$ on $C^1(\Delta^n,M)$ is said to have 
\emph{Lipschitz determination} if there is $L>0$ such that $\mu$ 
is determined on the subset of $C^1$-simplices whose Lipschitz
constant is smaller than~$L$. 
\item Let $\Cmealip_\ast(M)$ denote the set of signed Borel measures
on $C^1(\Delta^n,M)$ that have finite total variation and Lipschitz
determination.  Then $(\Cmealip_n(M))_{n\ge 0}$ forms a chain complex
whose elements are called \emph{Lipschitz measure chains}. The
differential is given by the alternating sum of push-forwards induced
by face maps
\citelist{\cite{zastrow}*{Corollary~2.9}\cite{ratcliffe}*{p.~539}}.
The total variation defines a norm on each of these chain groups.
\item The homology groups of $\Cmealip_\ast(M)$ are denoted 
by $\Hmealip_\ast(M)$. They are equipped with the quotient semi-norm. 
\end{enumerate}
\end{definition}

The Lipschitz determination condition ensures that the
function~$\sigma \mapsto \int \sigma^* \dvol_M$ is bounded on the
supports of the measure chains in question. Therefore, Lipschitz
measure chains can be evaluated against the volume form:

\begin{remark}\label{rem:eval is bounded Borel}
Let $M$ be a Riemannian $n$-manifold and let
$\mu\in\Cmealip_{n}(M)$. Then the function
\[I:C^1(\Delta^n,M)\rightarrow\bbR,~
\sigma\mapsto\langle\dvol_M,\sigma\rangle=\int_{\Delta^n}\sigma^\ast\dvol_M\]
is well defined, measurable, and $\mu$-almost everywhere bounded, thus $\mu$-integrable. 
We denote the integral $\int I d\mu$ by $\langle\dvol_M, \mu\rangle$. 
\end{remark}

\begin{definition}
For a Riemannian manifold $M$, we define the following subcomplex 
of~$\Clflip_*(M)$ 
(see Definition~\ref{def:locally finite lipschitz homology})
\begin{equation*}
\Clonelip_\ast(M)
  = \biggl\{ \sum_{i\in\bbN} a_i\sigma_i\in\Clflip_\ast(M)
           ;~\text{$\sigma_i$ smooth for all $i\in \bbN$, and}
             \sum_{i\in \bbN}\abs{a_i}<\infty
    \biggr\}.
\end{equation*}
A cycle in $\Clonelip_{\dim M}(M)$ is called a \emph{fundamental cycle} 
if it is a locally finite fundamental cycle in $\Clf_{\dim M}(M)$.  
\end{definition}

From now on, we refer to the setting in Setup~\ref{setup:smearing map}. 
Thurston's smearing technique is a cunning way of averaging
the simplices over the isometry group of the universal cover: 

\begin{proposition}\label{prop:smearing}
Let $\sigma:\Delta^i\rightarrow M$ be a smooth simplex, and let
$\ucov{\sigma}:\Delta^i\rightarrow U$ be a lift of~$\sigma$ to $U$. The
push-forward of $\mu_{\Lambda\bs G}$ under the map 
\begin{equation*}
\smear_{\ucov{\sigma}}:\Lambda\bs G\rightarrow C^1(\Delta^i,N),~\Lambda
g\mapsto p_N\circ g\ucov{\sigma}
\end{equation*}
does not depend on the choice of the lift of $\sigma$ and is denoted
by $\mu_\sigma$. Further there is a well-defined chain map 
\begin{equation*}
\smear_\ast:\Clonelip_\ast(M)\longrightarrow\Cmealip_\ast(N),
~\sum_\sigma a_\sigma\sigma\mapsto\sum_\sigma a_\sigma\mu_\sigma.
\end{equation*}
\end{proposition}

\begin{proof}
  One uses the 
  right $G$-invariance of~$\mu_{\Lambda\bs G}$ for 
  showing that $\smear_\ast$ is independent
  of the choice of the lifts and compatible with the boundary. 
  The computations are similar
  to the ones in the classical case~\cite{loehdiplom}*{Section~5.4}.
\end{proof}

In the proof of the proportionality principle
(Theorem~\ref{thm:proportionality most general}), it is essential to
be able to determine the map induced by smearing in the top homology. 
We achieve this by evaluating with respect to the volume form.

\begin{lemma}\label{lem:eval smearing at volume form}
For every fundamental cycle
$c\in\Clonelip_n(M)$ we have 
\begin{equation*}
\bigl\langle\dvol_N,\smear_n(c)\bigr\rangle=\int_{C^1(\Delta^n,N)}\int_{\Delta^n}\sigma^\ast\dvol_Nd\smear_n(c)(
\sigma)=\vol(M).
\end{equation*}
\end{lemma}

\begin{remark}\label{rem:smooth simplices for simplicial volume}
There exists a fundamental cycle in
$\Clonelip_n(M)$ if and only if
$\lipvol{M}<\infty$. 
Equivalently, the Lipschitz simplicial
  volume can be computed by smooth cycles:
  \begin{equation}\label{eq:by smooth cycles}
        \lipvol{M}
     =  \inf \bigl\{ \abs{c}_1
                   ;~\text{%
	             $c \in \Clf_n(M)$ smooth fundamental cycle, 
	             $\lipschitz(c) < \infty$}
             \bigr\}. 
  \end{equation}
This can be shown without curvature conditions 
using relative approximation theorems for Lipschitz maps by 
smooth ones but in the case of non-positively curved manifolds 
the straightening technique gives 
a quick proof of~\ref{eq:by smooth cycles}: 

If $c = \sum_{i \in I} a_i\sigma_i 
\in \Clf_*(M)$ satisfies~$\lipschitz(c) < \infty$, then
  Proposition~\ref{prop:geodesic simplex - non-positive curvature} 
  and Remark~\ref{rem:diameter of geodesic simplices}
  show that also the straightened chain~$c'=\sum_{i\in I}a_i \cdot
  s_M(\sigma_i)$ is both Lipschitz and locally finite. Moreover, 
it is smooth 
by~\ref{prop:geodesic simplex - non-positive curvature}. 
Thus, straightening chains gives rise 
to a chain map~$\Clflip_*(M) \rightarrow\Clflip_*(M)$. 
The same arguments as in 
Proposition~\ref{prop:geodesic straightening} and 
Lemma~\ref{lem:lemma from lee} also apply in 
the locally finite case, which implies that this chain 
map is homotopic to the identity. Hence 
$[c']=[c]$, which, combined with $\abs{c'}_1\le\abs{c}_1$, 
shows~\ref{eq:by smooth cycles}. 
\end{remark}

\begin{proof}[Proof of Lemma~\ref{lem:eval smearing at volume form}]
In view of Remark~\ref{rem:eval is bounded Borel}, the double integral
in the lemma is well-defined. Because the universal covering maps~$p_M$
and~$p_N$ are locally isometric, we obtain (where we write $c =
\sum_{\sigma} a_\sigma \sigma$)
\begin{align*}
\langle\dvol_N,\smear_n(c)\rangle
&=\sum_\sigma a_\sigma\langle\dvol_N, \mu_{\ucov{\sigma}}\rangle\\
&=\sum_\sigma
a_\sigma\int_{C^1(\Delta^n,N)}\langle\dvol_N, \rho\rangle \;d\mu_{\ucov{\sigma}}(\rho)\\
&=\sum_\sigma a_\sigma\int_{\Lambda\backslash G}\langle\dvol_N, p_N\circ
g\ucov{\sigma}\rangle \;d\mu_{\Lambda\backslash G}(g)\\
&=\sum_\sigma a_\sigma\int_{\Lambda\backslash G}\langle
\dvol_U,g\ucov{\sigma}\rangle \;d\mu_{\Lambda\backslash G}(g)\\ 
&=\sum_\sigma a_\sigma\int_{\Lambda\backslash G}\langle
\dvol_U,\ucov{\sigma}\rangle \;d\mu_{\Lambda\backslash G}(g)\\ 
&=\sum_\sigma a_\sigma\int_{\Lambda\backslash G}\langle
\dvol_M,\sigma\rangle \;d\mu_{\Lambda\backslash G}(g). 
\end{align*}
By Proposition~\ref{prop:evaluation at the fundamental class}, the last
expression equals~$\vol(M)$. 
\end{proof}

\subsection{Proof of Theorem~\ref{thm:proportionality most general}}
\label{subsec:general proportionality}
In order to prove the proportionality principle
(Theorem~\ref{thm:proportionality most general}), we proceed in the
following steps: 
\begin{enumerate}
  \item First we construct a $\Lambda$-equivariant partition of~$U$
        into Borel sets of small diameter and a corresponding
        $\Lambda$-equivariant $1$-net. 
  \item Using the $1$-net and a
        straightening procedure, we develop a discrete version of the
        smearing map -- i.e., a mechanism turning fundamental cycles
        on~$M$ into cycles on~$N$. This has some similarity with 
        the construction by Benedetti and Petronio~\cite{bp}*{p.~114f}. 
  \item By comparing the discrete smearing with the original
        smearing, integration enables us to identify which class the
        smeared cycle represents. 
  \item In the final step, we compute the $\ell^1$-norm of the smeared
        cycle, thereby proving the theorem. 
\end{enumerate}

\begin{proof}[Proof of Theorem~\ref{thm:proportionality most general}]
Like in the previous paragraphs, we refer to the notation established
in Setup~\ref{setup:smearing map}. 

\subsubsection{Construction of a suitable $\Lambda$-equivariant
  partition of~$U$ into Borel sets}
  
By locally subdividing a triangulation of~$N$, it is possible to
construct a locally finite (and hence countable) set~$T \subset N$ and
a partition~$(F_x)_{x \in T}$ of~$N$ into Borel sets with the
following properties: For each~$x\in T$ we have $x\in F_x$, the
diameter of~$F_x$ is at most~$1/2$ (thus, $T$ is a $1$-net in~$N$),
and the universal cover~$p_N$ is trivial over~$F_x$. 

Let $\ucov T \subset U$ be a lift of~$T$ to~$U = \ucov N$. In view of
the triviality condition, we find a corresponding
$\Lambda$-equivariant partition~$\ucov F := (\ucov F_x)_{x \in \Lambda \cdot
\ucov T}$ of~$U$ into Borel sets of diameter at most~$1/2$. Note 
that $\Lambda\cdot\ucov T$ is locally finite since $\Lambda$ acts 
properly on $\ucov N$.

\subsubsection{Discrete version of the smearing map}

In order to construct the discrete version of the smearing map, we
first define a version~$\str$ of the geodesic straightening that turns
simplices in~$U$ into geodesic simplices with vertices
in~$\Lambda\cdot \ucov T$: For an $i$-simplex
$\rho:\Delta^i\rightarrow U$ we define the geodesic simplex
\[      \str_i(\rho)
    :=  [x_0, \dots, x_i],
\]
where $x_0, \dots, x_i \in \Lambda \cdot \ucov T$ are the elements
uniquely determined by the requirement that for all~$j \in \{0,\dots,
i\}$ the $j$-th vertex of~$\rho$ lies in~$\ucov F_{x_j}$. By
Proposition~\ref{prop:geodesic simplex - non-positive curvature}, the
simplex~$\str_i(\rho)$ is smooth.  Because the partition~$\ucov F$ is
$\Lambda$-equivariant, so is~$\str_i$. Using the fact that all
elements of~$\ucov F$ are Borel and that $\Lambda \cdot \ucov T$ is
countable, it is not difficult to see that the map~$\str_i \colon
C^1(\Delta^i, U) \rightarrow C^1(\Delta^i,N)$ is Borel with respect to
the $C^1$-topology.  Moreover, for all~$k \in \{0,\dots,i \}$
\begin{equation}\label{eq:face compatibility}
\str_{i-1}(\partial_k\rho)=\partial_k\str_i(\rho).
\end{equation}

For $i\in \bbN$ we write  
\begin{align*}
       S_i 
  := &\,\bigl\{ p_N\circ\sigma
              ;~\sigma:\Delta^i\rightarrow U
                \text{ geodesic simplex with vertices in~$\Lambda\cdot\ucov T$}
        \bigr\}\\
  \subset &\,C^1(\Delta^i, N)
  ,
\end{align*}
and for every simplex $\sigma:\Delta^i\rightarrow U$ we define a map
\[f_\sigma: G\rightarrow S_i,~g\mapsto
p_N\circ\str_i(g\sigma);\] 
The map~$f_\sigma$ is Borel because $\str_i$ is Borel and the action
of~$G$ is $C^1$-continuous (the compact-open topology on $G$ coincides 
with the $C^1$-topology~\cite{loehdiplom}*{Theorem~5.12}). 
Furthermore, $f_\sigma$ induces a
well-defined Borel map $f_\sigma:\Lambda\backslash G\rightarrow
S_i$, which we denote by the same symbol.

We now consider
the following discrete approximation of the smearing map defined in
Proposition~\ref{prop:smearing} 
\begin{gather}
\phi_\ast:\Clonelip_\ast(M)\rightarrow\Clonelip_\ast(N)\notag\\ \label{eq:discretized smearing map}
\phi_i\biggl(\sum_{k\in \bbN}a_k\sigma_k\biggr):=\sum_{\rho\in S_i}\biggl(\sum_{k\in\bbN}a_k\cdot\mu_{\Lambda\bs
  G}\bigl(f_{\ucov{\sigma}_k}^{-1}(\rho)\bigr)\biggl)\cdot \rho
\end{gather}
where each $\ucov{\sigma}_k$ is a lift of $\sigma_k$ to $U$. First we
show that $\phi_\ast$ is well-defined: The number~$\mu_{\Lambda\bs
G}(f^{-1}_{\ucov{\sigma}}(\rho))$ does not depend on the choice of the
lift $\ucov{\sigma}$ of the simplex~$\sigma$ because $\mu_{\Lambda\bs
G}$ is invariant under right multiplication of $G$.  If
$L=\lipschitz(\sigma)$, any lift~$\ucov{\sigma}$ has diameter at most
$\sqrt 2 L$. Hence, each pair of vertices of $\str_i(g\ucov{\sigma})$
has distance at most~$1+\sqrt 2 L$. In view of
Proposition~\ref{prop:geodesic simplex - non-positive curvature} and
Remark~\ref{rem:diameter of geodesic simplices},
$\str_i(g\ucov{\sigma})$, and thus $f_{\ucov{\sigma}}(g)$, are 
smooth and have a Lipschitz constant
depending only on~$L$. Hence 
there is a uniform bound on the Lipschitz constants of simplices
appearing in the right hand sum of~\ref{eq:discretized smearing
map}. This also implies that \eqref{eq:discretized smearing map}
defines a locally finite chain because both $\Lambda \cdot \ucov T$
and $T$ are locally finite.  Therefore, $\phi_i$ is a well-defined
homomorphism for every $i\in \bbN$.

Next we prove that $\phi_\ast$ is a chain homomorphism:  
From~\ref{eq:face compatibility} we obtain
\begin{equation*}
\bigcup_{\rho\text{ with
  }\partial_k\rho=\xi}\bigl\{\Lambda g\in\Lambda\backslash
G;~p_N\circ\str_i(g\ucov{\sigma})=\rho\bigr\}
=\bigl\{\Lambda g\in\Lambda\backslash
G;~p_N\circ\str_{i-1}(g\partial_k\ucov{\sigma})=\xi\bigr\}
\end{equation*}
for all~$\sigma \in \map(\Delta^i,N)$, $k\in\{0,\dots,i\}$, and
all~$\xi \in\map(\Delta^{i-1},N)$.  Because the left hand side is a
disjoint, at most countable, union this implies that
\[\sum_{\rho\text{ with
  }\partial_k\rho=\xi}\mu_{\Lambda\bs
  G}\bigl(f_{\ucov{\sigma}}^{-1}(\rho)\bigr)=\mu_{\Lambda\bs
  G}\bigl(f_{\partial_k\ucov{\sigma}}^{-1}(\xi)\bigr).\]
Therefore, we deduce
\begin{align*}
\partial_k\phi_i(\sigma) &=\sum_{\rho\in S_i}\mu_{\Lambda\bs G}\bigl(f_{\ucov{\sigma}}^{-1}(\rho)\bigr)\cdot\partial_k\rho\\ 
&=\sum_{\xi\in S_{i-1}}\sum_{\rho\text{ with
  }\partial_k\rho=\xi}\mu_{\Lambda\bs G}\bigl(f_{\ucov{\sigma}}^{-1}(\rho)\bigr)\cdot\xi\\
&=\sum_{\xi\in
  S_{i-1}}\mu_{\Lambda\bs G}\bigl(f_{\partial_k\ucov{\sigma}}^{-1}(\xi)\bigr)\cdot\xi\\
&=\phi_{i-1}(\partial_k\sigma),
\end{align*}
which shows that $\varphi_*$ is a chain map.

\subsubsection{Comparison with the original smearing map}
Let 
\[ j_\ast : \Clonelip_\ast(N)\rightarrow \Cmealip_\ast(N)
\] 
be the chain map that is the obvious extension of the map given by
mapping a simplex~$\sigma$ to the atomic measure concentrated
in~$\{\sigma\}$.  Next we show that there is a chain homotopy between
the smearing map $\smear_\ast$ given in
Proposition~\ref{prop:smearing} and the
composition~$j_\ast\circ\phi_\ast$: For any smooth simplex
$\sigma:\Delta^i\rightarrow U$ and~$g\in G$ the geodesic homotopy from
$\str_i(g\sigma)$ to $g\sigma$ followed by $p_N$ defines a 
map~$h_\sigma(g):\Delta^i\times I\rightarrow N$. By
Proposition~\ref{prop:smooth join}, Proposition~\ref{prop:geodesic
simplex - non-positive curvature}, and Remark~\ref{rem:diameter of
geodesic simplices}, $h_\sigma(g)$ is smooth and 
its  Lipschitz constant is
bounded from above in terms of the Lipschitz constant of
$\sigma$. Moreover, Proposition~\ref{prop:smooth join} shows that the
map~$h_\sigma \colon G \rightarrow C^1(\Delta^i,N)$ is Borel with
respect to the $C^1$-topology.  Because $\str_\ast$ is
$\Lambda$-equivariant, we obtain a well-defined Borel map
\begin{equation*}
h_\sigma:\Lambda\backslash
G\rightarrow C^1(\Delta^i\times I,N) 
\end{equation*}
satisfying
\begin{align}\label{eq:homotopy}
h_\sigma(\Lambda g)\vert_{\Delta^i\times\{0\}}&=f_\sigma(g),\\
h_\sigma(\Lambda g)\vert_{\Delta^i\times\{1\}}&=p_N\circ g\sigma\nonumber. 
\end{align}

It is also clear that for each face map
$\partial_k:\Delta^{i-1}\rightarrow\Delta^i$ and every simplex
$\sigma:\Delta^i\rightarrow U$ we have
\[h_{\sigma\circ\partial_k}(\Lambda g)=h_\sigma(\Lambda g)\circ\bigl(\partial_k\times\id_I\bigr).\]
Retaining the notation of Lemma~\ref{lem:lemma from lee} and
Remark~\ref{rem:notation from lee}, 
for every $\sigma:\Delta^i\rightarrow U$ and every~$k \in \{0,\dots,i\}$ let 
$\nu_{\sigma,k}$ be the push-forward of $\mu_{\Lambda\bs G}$ 
under the map 
\[\Lambda\backslash G\rightarrow C^1(\Delta^{i+1},N),~\Lambda
g\mapsto h_\sigma(g)\circ G_{i,k}.\]
If $\sigma$ is a simplex in $M$ and $\ucov{\sigma}$ a lift to $U$, then 
$\nu_{\ucov{\sigma},k}$ does not depend on the choice of the lift
and will be also denoted by $\nu_{\sigma,k}$. 
We now define the homomorphism 
\begin{gather*}
H_\ast:\Clonelip_\ast(M)\rightarrow
\Cmealip_{\ast+1}(N),~H_i(\sigma):=\sum_{k=0}\nu_{\sigma,k}. 
\end{gather*}
Lemma~\ref{lem:lemma from lee} and~\ref{eq:homotopy} yield~\cite{zastrow}*{Theorem~2.1 (1)}
\begin{equation*}
\partial H_i(\sigma)+H_{i-1}\partial\sigma =j_i(\phi_i(\sigma))-\smear_i(\sigma)
\end{equation*}
for every $i$-simplex $\sigma$ in $M$. 
Thus $H_\ast$ 
is the desired chain homotopy
$j_\ast\circ\phi_\ast\simeq\smear_\ast$. 

The evaluation with $\dvol_N$ 
(cf.~Remark~\ref{rem:eval is bounded Borel}) 
is compatible with $j_\ast$, that is, 
\begin{equation*}
\bigl\langle \dvol_N, j_\ast(c)\bigr\rangle=\langle \dvol_N, c\rangle
\end{equation*}
for every $c\in\Clonelip_\ast(N)$. 

Let $c\in\Clonelip_n(M)$ be a fundamental cycle. Because evaluation
with $\dvol_N$ is well-defined on homology classes and by
Lemma~\ref{lem:eval smearing at volume form}, we obtain that
\begin{align*}
\bigl\langle\dvol_N, \phi_n(c)\bigr\rangle&=\bigl\langle
\dvol_N,j_n(\phi_n(c))\bigr\rangle\\
&=\bigl\langle\dvol_N,\smear_n(c)\bigr\rangle\\
&=\vol(M).
\end{align*}
Now Proposition~\ref{prop:evaluation at the fundamental class} lets 
us determine the homology class of $\phi_n(c)$ as 
\begin{equation}\label{eq:detect fundamental class}
[\phi_n(c)]=\frac{\vol(M)}{\vol(N)}\cdot[N]. 
\end{equation}

\subsubsection{The norm estimate and conclusion of proof}
By symmetry we only have to show that 
\[\frac{\lipvol{M}}{\vol(M)}\ge\frac{\lipvol{N}}{\vol(N)},\]
and in addition we can assume $\lipvol{M}<\infty$. By
Remark~\ref{rem:smooth simplices for simplicial volume}, 
we can compute the Lipschitz simplicial volume~$\lipvol M$
by fundamental cycles lying in the chain complex~$\Clonelip_*(M)$.
Let $c=\sum_{k\in \bbN}a_k\sigma_k\in\Clonelip_n(M)$ be a fundamental
cycle of~$M$. 
Because of~\ref{eq:detect fundamental class} it suffices to show that 
\begin{equation*}
\abs{\phi_n(c)}_1\le\abs{c}_1,
\end{equation*}
which is a consequence the following computation: 
\begin{align*}
\abs{\phi_n(c)}_1&\le\sum_{\rho\in S_n}\sum_{k\in\bbN}\abs{a_k}\cdot\mu\bigl(f_{\ucov{\sigma_k}}^{-1}(\rho)\bigr)\\
&=\sum_{k\in \bbN}\sum_{\rho\in S_n}\abs{a_k}\cdot\mu\bigl(f_{\ucov{\sigma_k}}^{-1}(\rho)\bigr)\\
&=\sum_{k\in \bbN}\abs{a_k}\\
&=\abs{c}_1.
\end{align*}
This finishes the proof of the proportionality principle. 
\end{proof}

\section{Vanishing results for the locally finite simplicial volume}
\label{sec:vanishing results}

\noindent
In this section, we give a proof of the vanishing theorem
(Theorem~\ref{thm:vanishing of locally finite simplicial volume}); the
proof is based on the fact that locally symmetric spaces of higher
$\bbQ$-rank admit ``amenable'' coverings of sufficiently small
multiplicity and Gromov's vanishing finiteness theorem.

As a first step, we recall Gromov's definition of amenable subsets and
sequences of subsets that are amenable at
infinity~\cite{gromov}*{p.~58} and his vanishing-finiteness theorem:

\begin{definition} Let $X$ be a topological space. 
\begin{enumerate}
\item A subset $U\subset X$ is called \emph{amenable in $X$} 
if for every basepoint~$x \in U$ the
subgroup~$\im\bigl(\pi_1(U,x)\rightarrow\pi_1(X,x)\bigr)$ is
amenable. 
\item A sequence $(U_i)_{i\in \bbN}$ of subsets of~$X$ is called 
$\emph{amenable at infinity}$ if there is an increasing sequence 
of compact subsets $(K_i)_{i\in \bbN}$ of~$X$ with $U_i\subset X-K_i$, 
$X=\bigcup_{i \in \bbN} K_i$, and such that $U_i$ is amenable in~$X-K_i$ for
sufficiently large~$i\in\bbN$.  
\end{enumerate}
\end{definition}

\begin{theorem}[Vanishing-finiteness theorem for simplicial
  volume~\cite{gromov}*{Corollary~(A) on p.~58}]\label{thm:vanishing-finiteness} 
Let $M$ be a manifold without boundary of dimension~$n$. 
Let $(U_i)_{i\in\bbN}$ be a locally finite covering of~$M$ by open, relatively 
compact subsets such that each point of~$M$ is contained in at most~$n$ such subsets. 
If every $U_i$ is amenable in~$M$ and $(U_i)_{i\in \bbN}$ is amenable at infinity, 
then $\lfvol{M}=0$. 
\end{theorem}

As a next step, we provide a construction of locally finite coverings
with small multiplicity by relatively compact, open, amenable subsets;
notice however that such a covering is not necessarily amenable at
infinity. 

\begin{theorem}\label{thm:construction of amenable covers}
  Let $M$ be a manifold and $\Gamma=\pi_1(M)$.  Assume that $\Gamma$
  admits a finite model for its classifying space $B\Gamma$ of
  dimension $k$.  Then there is a locally finite covering of $M$ by
  relatively compact, amenable, open subsets such that every point of
  $M$ is contained in at most $k+2$ such subsets.
\end{theorem}

\begin{proof}
Since $B\Gamma$ is $k$-dimensional and compact, 
every open covering of~$B\Gamma$ has a
finite refinement with multiplicity at 
most~$k+1$~\cite{hure}*{Theorem~V~1~on p.~54}. 
Starting with a covering of~$B\Gamma$ by open, contractible sets, let
$(V_j)_{j\in J}$ be a finite refinement of multiplicity at
most~$k+1$. 

We pull this covering back to~$M$ via the classifying map~$\phi
\colon M \rightarrow B\Gamma$: For~$j \in J$ let
\[ U_j := \phi^{-1}(V_j).
\] 
By construction, $(U_j)_{j \in J}$ is an open covering of~$M$ with
multiplicity at most~$k + 1$. However, the sets~$U_ j$ may not be
relatively compact.

To achieve a nice covering of~$M$ by relatively compact sets, we
combine the covering~$(U_j)_{j \in J}$ with another covering of~$M$ of
small multiplicity consisting of relatively compact sets, which is
constructed as follows: For every~$j \in J$ we choose a
covering~$\calr_j$ of~$\bbR$ by bounded, open intervals such that
each~$\calr_j$ has multiplicity~$2$ and for~$i \neq j$ the
cover~$\calr_i \smallcoprod \calr_j$ (disjoint union) has multiplicity at
most~$3$. This is possible because $J$ is finite.

Let $f \colon M \rightarrow \bbR$ be a proper function. 
We show now that the combined covering 
\[\calu := \bigl( U_j\cap f^{-1}(W)
           \bigr)_{j\in J,~W\in\calr_j}
\] 
of~$M$ has the desired properties: In the following, by definition, we
say that the \emph{J-index} of~$U_j\cap f^{-1}(W)$ is~$j$.

Because $f$ is proper and the elements of the~$\calr_j$ are bounded,
each set in~$\calu$ is relatively compact.
 
Since $\phi \colon \pi_1(M) \rightarrow \pi_1(B\Gamma)$ is an
isomorphism, the inclusion~$U_j\cap f^{-1}(W)\hookrightarrow
M$ is trivial on the level of~$\pi_1$ if and only if its composition
with $\phi$ is so. But the composition with~$\phi$ factors over the
inclusion $V_j\hookrightarrow B\Gamma$, which is trivial
in~$\pi_1$. In particular, each element of $\calu$ is an amenable subset
of~$M$.

It remains to verify that~$\calu$ has multiplicity at most~$k + 2$:
Suppose there is a subset~$\calu_0\subset\calu$ of $k+3$ sets whose
intersection is non-empty. Because the elements of~$\calu_0$ have at
most $k + 1$ different $J$-indices, and the multiplicity of each of
the~$\calr_j$ is at most~$2$, there must be $i \neq j \in J$ such that
there are at least two elements in~$\calu_0$ having $J$-index~$i$, and
at least two with $J$-index~$j$. But this contradicts the fact that
$\calr_i\smallcoprod\calr_j$ has multiplicity at most~$3$. So the
multiplicity of~$\calu$ is at most~$k+2$.
\end{proof}

In order to obtain a suitable amenable covering that is amenable at
infinity, we impose additional constraints on the fundamental group of
the boundary;
one should compare this also with Gromov's remark on
subpolyhedra~\cite{gromov}*{p.~59}.

\begin{corollary}\label{cor: small classifying space}
	Let $M$ be the interior of a compact, $n$-dimensional manifold $W$ with 
	boundary $\partial W$. Assume that $B\pi_1(M)$ 
  	admits a finite model of dimension at most~$n-2$ and that  
	at least one of the following conditions is satisfied 
	\begin{enumerate}
		\item The fundamental group~$\pi_1(\partial W; x)$ is amenable for all $x\in\partial W$. 
		\item For all~$x \in \partial W$ the inclusion induces
                  an injection~ $\pi_1(\partial W;x)\rightarrow
                  \pi_1(W;x)$.  
	\end{enumerate}
	Then $\lfvol{M}=0$.
\end{corollary}
\begin{proof}
  	By Theorem~\ref{thm:construction of amenable covers} 
	we obtain a covering $(U_i)_{i\in \bbN}$ of~$M$ by open, relatively compact,  
	amenable subsets in~$M$ which has multiplicity~$\le n$. 
	Let $(V_i)_{i\in \bbN}$ be a decreasing 
	sequence of open neighborhoods in~$W$ of the
        boundary~$\partial W$ with~$\bigcap_{i \in \bbN} V_i =
        \partial W$ and $\bigcup_{i \in \bbN} V_i = W$.  By choosing
        collar neighborhoods of~$\partial W$ we can assume that
        $\partial W$ is a deformation retract of~$V_i$ for all
        large~$i \in \bbN$. Because $(U_i)_{i \in \bbN}$ is locally
        finite, we additionally can assume that $U_i \subset V_i$ for
        all~$i \in \bbN$.
	
	If $\pi_1(\partial W;x)$ is amenable for every basepoint $x$
        then $U_i$ is obviously an amenable subset of~$V_i$ for all
        large~$i \in \bbN$. 
	If the inclusion maps $\partial W\rightarrow W$ are
        $\pi_1$-injective then so are the inclusion maps $V_i\cap
        M\rightarrow M$ for all large~$i \in \bbN$, and the
        amenability of the subset $U_i\subset V_i\cap M$ follows from
        the one of $U_i\subset M$.
	
	In either case we can now apply Gromov's vanishing-finiteness 
	theorem~\ref{thm:vanishing-finiteness}. 
\end{proof}

\begin{example}[Products of open manifolds with non-zero simplicial
  volume]\label{ex:products with non-zero sv}
  We consider the open manifold~$M := W^\circ \times \bbR$, where 
  $(W,\partial W)$ is the surface with boundary obtained by removing a
  finite number of pairwise disjoint, open discs from an oriented, closed,
  connected surface of genus at least~$1$. 

  Then the fundamental group~$\pi_1(M) \cong \pi_1(W)$ is a finitely
  generated free group and thus admits a finite model of dimension~$1
  = \dim M -2$. 

  However, we can view~$M$ as the interior of the compact
  manifold~$W \times [0,1]$ whose boundary is nothing but an oriented,
  closed, connected surface of genus at least~$2$; in particular, this
  boundary has non-zero simplicial volume, which forces the simplicial
  volume of~$M$ to be
  infinite~\citelist{\cite{gromov}*{p.~17}\cite{loehphd}*{Corollary~6.2}}.

  In fact, tracking down the construction of an open covering in the
  proof of Theorem~\ref{thm:construction of amenable covers} shows
  that this particular covering is amenable but not amenable at
  infinity.
\end{example}

In particular, the finiteness hypothesis in the corollary is not
sufficient for the vanishing of the simplicial volume. The following
cohomological criterion helps to check whether the finiteness
hypothesis in the corollary is satisfied.

\begin{lemma}\label{lem:brown stuff}
  Let $\Gamma$ be a group that has a finite
  model for its classifying space $B\Gamma$. 
  \begin{enumerate}
    \item If $\cd \Gamma \neq 2$, then there is a finite model
          of~$B\Gamma$ whose dimension equals the integral
          cohomological dimension~$\cd \Gamma$ of $\Gamma$. 
    \item If $\cd \Gamma = 2$, then there is a finite model 
          of~$B \Gamma$ of dimension at most~$3$.
  \end{enumerate}
\end{lemma}

\begin{proof}
  Because there is a finite model for~$B\Gamma$, the group~$\Gamma$ is
  finitely presented and of type~FL. Therefore, a classic
  result of Eilenberg and Ganea shows that there is a finite model
  of~$B\Gamma$ of dimension~$\max\{\cd \Gamma,
  3\}$~\cite{brown}*{Theorem~VIII.7.1}.

  If $\cd \Gamma = 0$, then $\Gamma$ is the trivial group and hence
  the one-point space is a model for~$B \Gamma$. If $\cd \Gamma = 1$,
  then $\Gamma$ is free by a theorem of Stallings and
  Swan~\cite{stallings}; because $\Gamma$ is finitely presented, $\Gamma$
  is a finitely generated free group. In particular, we can take a
  finite wedge of circles as a finite, one-dimensional model
  for~$B\Gamma$.
\end{proof}

Using the techniques established in this section, we prove the
vanishing theorem for the locally finite simplicial volume of
non-compact locally symmetric spaces of $\bbQ$-rank $\ge 3$ 
(Theorem~\ref{thm:vanishing of locally finite simplicial volume}):

\begin{proof}[Proof of Theorem~\ref{thm:vanishing of locally finite simplicial volume}] 
  The locally symmetric space $M=\Gamma\bs X$ is a model of $B\Gamma$
  because $X$ is non-positively curved~\cite{eberlein}*{Sections~2.1 and~2.2}, 
  thus contractible. 
  Moreover, $M$ is homotopy equivalent to the Borel-Serre
  compactification~$W$ of~$M$~\cite{borelserre}, which thus is a
  finite model of $B\Gamma$. For $\rk_\bbQ\Gamma\ge 3$ 
  the inclusion $\partial W\rightarrow W$ is a 
  $\pi_1$-isomorphism~\cite{block+weinberger}*{Proposition~2.3}. 
  Furthermore, we
  have~\cite{borelserre}*{Corollary~11.4.3}
  \begin{align*}\label{eq:cd and dim} 
    \cd \Gamma  =    \dim X - \rk_\bbQ \Gamma. 
  \end{align*}
  Therefore, Lemma~\ref{lem:brown stuff} shows that there is a finite
  model for~$B \Gamma$ of dimension at most~$\max\{\dim X-\rk_\bbQ
  \Gamma, 3\}\le\dim X-2$.
  
  Thus, Corollary~\ref{cor: small classifying space} yields 
  the vanishing of~$\lfvol M$. 
\end{proof}

\end{document}